\documentclass[12pt]{article}

\usepackage{amsmath}
\usepackage{amsfonts}
\usepackage{latexsym}
\usepackage{graphicx}
\usepackage{amsthm}

\newtheorem{theorem}{Theorem}[section]

\newtheorem{proposition}[theorem]{Proposition}

\newtheorem{remark}[theorem]{Remark}

\date{ }

\begin{document}

\title{Topology optimization and boundary observation
for clamped plates}

\author{Cornel Marius Murea$^1$, Dan Tiba$^2$\\
{\normalsize $^1$ D\'epartement de Math\'ematiques, IRIMAS,}\\
{\normalsize Universit\'e de Haute Alsace, France,}\\
{\normalsize cornel.murea@uha.fr}\\
{\normalsize $^2$ Institute of Mathematics (Romanian Academy) and}\\ 
{\normalsize Academy of Romanian Scientists, Bucharest, Romania,}\\ 
{\normalsize dan.tiba@imar.ro}
}

\maketitle

\begin{abstract}
   We indicate a new approach to the optimization of the clamped plates with holes.
  It is based on the use of Hamiltonian systems and the penalization of the performance
  index. The  alternative technique employing the penalization of the state system, cannot
  be applied in this case due to the (two) Dirichlet boundary conditions. We also include numerical tests
exhibiting both shape and topological modifications, both creating and closing holes.

\vspace{3mm}

\vspace{2mm}
\textbf{MCS 2020:} 49Q10; 90C90
\end{abstract}

\section{Introduction}
\setcounter{equation}{0}

Geometric optimization problems appeared very early in the development
of mathematics and we quote just the famous Dido's problem, almost 3000 years
old, Bandle \cite{Bandle2017}, Alekseev, Tikhomirov, Fomin \cite{Alekseev1987}.
In the beginning of the previous century, Hadamard \cite{Hadamard1909}
introduced domain variations that were later generalized and play a fundamental role
in shape optimization \cite{Delfour2011}, \cite{Sokolowski1992}. The literature
devoted to this subject and to various applications is very rich and we quote just
the monograph \cite{Pironneau1984}, \cite{Henrot2005}, \cite{NS_Tiba2006},
\cite{Haslinger1996}, \cite{Haslinger2003} and their references.
Topology optimization was investigated mainly by the topological asymptotics
and the topological gradient methods in \cite{Amstutz2006}, \cite{Amstutz2017},
\cite{Masmoudi2007}, \cite{Wang2006} and in the monographs \cite{Plotnikov2012},
\cite{Novotny2013}, \cite{Novotny2019}.
Fixed domain methods like integrated topology optimization \cite{Bendsoe2003},
homogenization \cite{Allaire2007}, penalization \cite{NP_Tiba2009}, \cite{N_Tiba2012}
also allow topological optimization, in combination with shape optimization.
In particular, the paper \cite{MT2019a} of the authors discusses a penalization
approach to the optimization of a simply supported plate with holes, including
numerical experiments. Unfortunately, this methodology seems not possible
to be extended for clamped plates since the penalization of the state system seems
very much dependent on the
boundary conditions. See \cite{Zhou2017}, \cite{Zhou2014} where it is shown
that the penalization techniques for Dirichlet, respectively Neumann boundary
conditions, are very different and combining them (as necessary in the case of clamped plates) seems not possible.

Recently, a new local representation of the unknown geometries, called
the implicit parametrization approach, based on the use of certain Hamiltonian systems, was introduced in \cite{Tiba2013},
\cite{N_Tiba2015}, \cite{Tiba2018}.
In dimension two, which is a case of interest in optimal design and appropriate to the study of plates, the obtained
parametrization is even global \cite{Tiba2018a}, under certain assumptions.
New applications to topology optimization problems, including boundary observation, are possible \cite{MT2022}, \cite{MT2022a}.

In this work we show that the Hamilton approach, using implicit parametrizations,
provides a solution to the shape/topology optimization problems associated with
clamped plates. For other results concerning shape optimization problems
governed by plate models, we quote \cite{Sokolowski1992} sect 3.7, 
\cite{Novotny2013} sect 3.5 and 4.3, \cite{Delfour2011} sect 9.4.2.
Thickness optimization problems and other geometric optimization problems
expressed as control by the coefficients problems are investigated in 
\cite{Arnautu2000}, \cite{NS_Tiba2006} Ch. 6 (including curved structures like shells).

We underline that our techniques in this paper and in the previous works, don't need changes of
variables in the domains appearing in various iterations, and this allows combined shape
and topology variations.

The plan of the paper is as follows. In the next section we collect some
preliminaries and we formulate the topology optimization problem.
Section 3 gives general approximation results and discusses the applicability
of the gradient descent methods. Section 4 continues the analysis via the finite
element method and at the discrete level. The last section is devoted to some relevant
numerical experiments.

\section{Preliminaries and problem formulation}
\setcounter{equation}{0}

We denote by $\mathcal{O}$ the collection of the admissible domains, not
supposed to be simply connected. They are all contained in a given bounded
domain $D\subset \mathbb{R}^2$ and we introduce the clamped plate model:
\begin{eqnarray}
  \Delta h^3\Delta y_\Omega &=& f\hbox{ in }\Omega\in\mathcal{O},\label{2.1}\\
  y_\Omega=\frac{\partial y_\Omega}{\partial \mathbf{n}}&=& 0
  \hbox{ on }\partial\Omega,\label{2.2}
\end{eqnarray}  
where $h\in L^\infty(D)$ is the thickness of the plate, $h\geq m>0$ a.e. in $D$,
$f\in L^2(D)$ is the vertical load and $y_\Omega$ is the vertical deflection.
It is known that $y_\Omega\in H^2_0(\Omega)$ and can be obtained as a weak solution
of (\ref{2.1}), (\ref{2.2}), \cite{NS_Tiba2006}:
\begin{equation}\label{2.3}
  \int_\Omega h^3\Delta y_\Omega\Delta\varphi d\mathbf{x}=
  \int_\Omega f\varphi d\mathbf{x},\ \forall\varphi\in H^2_0(\Omega).
\end{equation}
Under smoothness assumptions on $h$ and $\Omega$ (class $\mathcal{C}^4$), then
the solution of (\ref{2.3}) satisfies $y_\Omega\in H^4(\Omega)\cap H^2_0(\Omega)$,
i.e. it is a strong solution \cite{Grisvard1985} Ch. 7, \cite{Agmon1959}. The same regularity
is imposed on $\partial D$ as well.

In geometric optimization problem, the functional variational approach
\cite{N_Tiba2012}, \cite{MT2022} assumes that the set of the admissible elements
$\mathcal{O}$ is defined via a family $\mathcal{F}\subset \mathcal{C}(\overline{D})$
of admissible controls:
\begin{equation}\label{2.4}
\Omega=\Omega_g=int\{ \mathbf{x}\in D;\ g(\mathbf{x})\leq 0 \},
\ g\in \mathcal{F}.
\end{equation}
Since $\Omega_g$ defined by (\ref{2.4}) is an open set, not necessarily connected,
we select the \textit{domain} $\Omega_g$ by the condition
\begin{equation}\label{2.5}
  \mathbf{x}_0\in\partial\Omega_g,\ \forall g\in \mathcal{F},
\end{equation}
where $\mathbf{x}_0\in D$ is some given point. Clearly, (\ref{2.4}), (\ref{2.5})
give $g(\mathbf{x}_0)=0$, for all $ g\in \mathcal{F}$, that has to be
imposed on $\mathcal{F}$, in order that the above procedure works.
These domains are just of class $\mathcal{C}$ \cite{Adams1975}, \cite{Tiba2013a}.
To obtain more regularity for $\partial\Omega_g$, we ask
$\mathcal{F}\subset \mathcal{C}^1(\overline{D})$ and
\begin{equation}\label{2.6}
| \nabla g(\mathbf{x}) | > 0,\ \forall g\in \mathcal{F},
\ \forall \mathbf{x}\in D, \hbox{with }g(\mathbf{x})=0.
\end{equation}
Then, $\partial\Omega_g$ is a class $\mathcal{C}^1$ by the implicit functions
theorem and
\begin{equation}\label{2.7}
  \partial\Omega_g=\{ \mathbf{x}\in D;\ g(\mathbf{x})= 0 \},
  \ \mathbf{x}_0\in \partial\Omega_g,\ \forall g\in \mathcal{F},
\end{equation}
and $\Omega_g$ is not necessarily simply connected.
More regularity for $\partial\Omega_g$ may be obtained if more regularity is
supposed for $\mathcal{F}$ and (\ref{2.6}) is valid. We also assume that
\begin{equation}\label{2.8}
g(\mathbf{x})>0,\ \forall\mathbf{x}\in \partial D,\ \forall g\in \mathcal{F}.
\end{equation}
Notice that (\ref{2.8}) simply ensures $\partial\Omega_g\cap \partial D=\emptyset$,
for any $g\in \mathcal{F}$. This framework was developed
mainly in \cite{MT2022}, \cite{Tiba2018a} and a central result is

\begin{theorem}\label{theo:2.1}
Under hypotheses (\ref{2.5}), (\ref{2.6}), (\ref{2.7}), (\ref{2.8}),
$\partial\Omega_g$ has a finite number of connected components, for
any $g\in \mathcal{F}$. Moreover, the connected component containing
$\mathbf{x}_0\in D\subset \mathbb{R}^2$, $g(\mathbf{x}_0)=0$, for all
$g\in \mathcal{F}$, is globally parametrized by the Hamiltonian system:
\begin{eqnarray}
z_1^\prime(t) & = & -\partial_2 g\left( z_1(t),z_2(t)\right),\ t\in I_g,\label{2.9}\\
z_2^\prime(t) & = &  \partial_1 g\left( z_1(t),z_2(t)\right),\ t\in I_g,\label{2.10}\\
\left( z_1(0),z_2(0)\right) & = & \mathbf{x}_0,\label{2.11}
\end{eqnarray}
that has a unique periodic solution and $I_g=[0,T_g]$ is the main period.
\end{theorem}

In general, the parametrization of manifolds, obtained by Hamiltonian-type systems,
has a local character \cite{Tiba2018}. Here, the system (\ref{2.9})-(\ref{2.10})
gives a global parametrization of $\partial\Omega_g$, due to the periodicity
of its solution, obtained via the Poincar\'e - Bendixson theory, \cite{Tiba2018a}.
Notice that the solution of Hamiltonian-type systems has the
uniqueness property, although the right-hand side is just continuous 
\cite{Tiba2018}. Global parametrizations of any component
of $\partial\Omega_g$ can be similarly obtained by fixing some initial condition
on it, in (\ref{2.9})-(\ref{2.11}). The corresponding period naturally
depends on each connected component.

It is possible to compute the directional derivative of $T_g$ with respect to
functional variations \cite{NP_Tiba2009}, \cite{N_Tiba2012}, $g+\lambda r$,
$g,r\in \mathcal{F}$ and
$\lambda \in \mathbb{R}$. We denote by $T_\lambda>0$ the main period of the
perturbed Hamilton system (\ref{2.9})-(\ref{2.11}) associated to $g+\lambda r$.
We have $T_\lambda\rightarrow T_g$ as $\lambda\rightarrow 0$ (see \cite{MT2019})
and

\begin{proposition}[\cite{MT2022}]\label{prop:2.2}
Under hypotheses (\ref{2.6})-(\ref{2.8}), we have
$$
\lim_{\lambda\rightarrow 0}
\frac{T_\lambda-T_g}{\lambda}=-\frac{w_2(T_g)}{z_2^\prime(T_g)}
$$
if $z_2^\prime(T_g)\neq 0$.
\end{proposition}
\noindent
Here $w=[w_1,w_2]$ is the solution of the system in variations associated to
(\ref{2.9})-(\ref{2.11}) via functional variations $g+\lambda r$, see
 Proposition 3.5.

We associate with the state system (\ref{2.1})-(\ref{2.2}) or, equivalently,
with its weak formulation (\ref{2.3}) the following general optimal design
problem:
\begin{equation}\label{2.12}
\min_{\Omega \in \mathcal{O}} J(y_\Omega).
\end{equation}

Since $\mathcal{O}$ is defined via (\ref{2.4}), starting from the admissible
family of functions $\mathcal{F}\subset \mathcal{C}^1(\overline{D})$ satisfying
(\ref{2.5})-(\ref{2.8}), relation (\ref{2.12}) may be rewritten as
\begin{equation}\label{2.13}
\min_{g \in \mathcal{F}} J(y_{\Omega_g}).
\end{equation}
subject to (\ref{2.1})-(\ref{2.2}).
In this way, the shape and topology optimization problem (\ref{2.12}),
(\ref{2.1}), (\ref{2.2}) will be transformed into an optimal control problem
defined in $D$. This will be performed in Section 3, together with a corresponding
approximation method, following techniques developed by the authors in
\cite{MT2022}. More constraints may be added to the above
optimization problem. For instance, if we want that all the admissible domains
contain some given open set $E\subset \Omega$, for all $\Omega\in \mathcal{O}$, then
we impose on $\mathcal{F}$ the condition
\begin{equation}\label{2.14}
g(\mathbf{x}) < 0,\ \forall \mathbf{x}\in E,\ \forall g \in \mathcal{F}.
\end{equation}
We also indicate here an example of boundary cost
functional (\ref{2.12}) or (\ref{2.13}) that shows the role played by
Thm. \ref{theo:2.1} (the assumptions will be detailed later):
\begin{eqnarray}
J(y_{\Omega}) & = &J(y_{\Omega_g})
  =\int_{\partial\Omega_g} j\left(\sigma,\Delta y_{\Omega_g}(\sigma)\right)d\sigma \nonumber\\
  &=& \int_0^{T_g} j\left(z_1(t),z_2(t),\Delta y_{\Omega_g}(z_1(t),z_2(t))\right)
  |(z_1^\prime(t),z_2^\prime(t)) | dt,
  \label{2.15}
\end{eqnarray}
where $z_1,z_2$ satisfy (\ref{2.9})-(\ref{2.11}). We show that $y_{\Omega_g}$
can be approximated
in $D$, as explained in the next section, and the main period $T_g$ can be computed
in a simple way, in the numerical
examples, see the last section. 

Notice that $\partial\Omega_g$ may have several connected components, according to
Thm. \ref{theo:2.1}. Then, the integrals in (\ref{2.15}) become sums of integrals,
each defined on another component. As already mentioned, on each component some
point (initial condition) has to be
fixed and an associated Hamiltonian system (\ref{2.9}) - (\ref{2.11}) has to be solved.
The solutions are all periodic due to (\ref{2.6}) and the corresponding periods may
be different.
Similar remarks are valid for the constraints (\ref{3.5}), (\ref{3.6}) and the
penalized cost (\ref{3.7}), in Section 3.

The control approach that we employ in this work
has a fixed domain character, that is all the computations are performed
in the given holdall bounded domain $D\subset\mathbb{R}^2$.

\section{The equivalent control problem and its approximation}
\setcounter{equation}{0}

We assume now that the thickness $h\equiv 1$ in $D$, for the sake of
simplicity, to avoid the regularity questions related to $h$,
\cite{Agmon1959}. Thickness optimization problems (for plates of a given shape) were
discussed in \cite{Arnautu2000}, \cite{NS_Tiba2006}. We underline that combined
thickness/shape/topology optimization problems would certainly be  of interest to
be investigated in some future paper.

We impose
$\mathcal{F}\subset \mathcal{C}^4(\overline{D})$ and (\ref{2.6}) in order that
$\partial\Omega_g$, defined in (\ref{2.7}) is of class $\mathcal{C}^4$, as required
in the classical paper \cite{Agmon1959}.
This ensures, for $f\in L^2(D)$, that the clamped plate model
(\ref{2.1}), (\ref{2.2}) has a strong solution
$y_g\in H^4(\Omega_g)\cap H^2_0(\Omega_g)$. The same regularity class $\mathcal{C}^4$ is assumed for $\partial D$. Denote by
$\tilde{y}_g\in H^4(D\setminus\overline{\Omega_g}) \cap H^2_0(D\setminus\Omega_g)$
some extension of $y_g\in D$, ensured by the trace theorem. Notice that
the concatenation of $y_g$, $\tilde{y}_g$, that we denote by $y_D$, satisfies
in $D$, the equation
\begin{eqnarray}
\Delta\Delta y_D & = & \chi_{\Omega_g} f + h_{\Omega_g} \hbox{ a.e. in }D,
\label{3.1}\\
y_D & = & \frac{\partial y_D}{\partial \mathbf{n}}=0\hbox{ on }\partial D,
\label{3.2}
\end{eqnarray}
where $\chi_{\Omega_g}$ is the characteristic function of $\Omega_g$,
$h_{\Omega_g}\in L^2(D)$ and $h_{\Omega_g}\equiv 0$ in $\Omega_g$,
$h_{\Omega_g}=\Delta\Delta \tilde{y}_g$ in $D\setminus\Omega_g$. 
Clearly, $\tilde{y}_{\Omega_g}$ and $h_{\Omega_g}$ are not unique. Due to 
(\ref{2.4}), (\ref{2.7}), by modifying $g\in\mathcal{F}$ outside of
$\Omega_g$, if necessary, we may assume that $g>0$ in
$D\setminus\overline{\Omega_g}$ and write (\ref{3.1}) in the form
\begin{equation}\label{3.3}
\Delta\Delta y_g = f + g_+^2 u \hbox{ a.e. in }D,
\end{equation}
where $g_+$ is the positive part of $g$ and has support equal with
$D\setminus\Omega_g$ and $u$ is  measurable in $D$,
$u=(\Delta\Delta \tilde{y}_g - f) /g_+^2$ in $D\setminus\overline{\Omega}_g$.
Consequently, $g_+^2 u\in L^2(D)$ and $y_g\in H^4(D)\cap H^2_0(D)$ by \cite{Agmon1959}. 
The state system is defined in $D$ by (\ref{3.2}), (\ref{3.3}) and we shall use
the notation $y_g\in H^4(D)\cap H^2_0(D)$ for the state, that is we neglect the
dependence on $u$ in (\ref{3.2}), (\ref{3.3}) and on $D$, that is fixed.

We concentrate on the cost functional (\ref{2.15}) which
corresponds to the difficult case of boundary observation and is specific to fourth order elliptic systems. Other cost functionals are discussed, for
instance, in \cite{Tiba2018a}, \cite{MT2022}, for second order problems
and can be easily adapted here.
Notice that $\Delta y_g \in H^2(D)\subset \mathcal{C}(\overline{D})$ by
the Sobolev embedding theorem \cite{Adams1975} and the cost functional
(\ref{2.15}) makes sense if $j$ is a Carath\'eodory function. We rewrite it
in the form
\begin{equation}\label{3.4}
\int_0^{T_g} j\left(z_1(t),z_2(t),\Delta y_g(z_1(t),z_2(t))\right)
  |(z_1^\prime(t),z_2^\prime(t)) | dt,
\end{equation}  
where $(z_1,z_2)\in \mathcal{C}^1(0,T_g)^2$ satisfies (\ref{2.9})-(\ref{2.11}) and
$y_g\in H^4(D)\cap H^2_0(D)$ is the solution of (\ref{3.2}), (\ref{3.3}).

The problem (\ref{3.2})-(\ref{3.4}), (\ref{2.9})-(\ref{2.11}) is an optimal
control problem defined in $D$, involving two controls $g\in \mathcal{F}$ and
$u$ that is measurable, such that $g_+^2 u\in L^2(D)$. We complete it with two constraints:
\begin{eqnarray}
&&
  \int_0^{T_g}
  \left| y_g\left(z_1(t),z_2(t)\right) \right|^2
  |(z_1^\prime(t),z_2^\prime(t)) | dt =0,
\label{3.5}\\
&&
  \int_0^{T_g}
  \left| \nabla y_g\left(z_1(t),z_2(t)\right)\cdot
    \frac{\nabla g\left(z_1(t),z_2(t)\right)}{|\nabla g\left(z_1(t),z_2(t)\right)|}
  \right|^2
  |(z_1^\prime(t),z_2^\prime(t)) | dt =0.
\label{3.6}
\end{eqnarray}
Notice that, due to (\ref{2.6}), conditions (\ref{3.5}), (\ref{3.6}) make
sense and they express the boundary conditions (\ref{2.2}) on $\partial\Omega_g$.
The constrained optimal control problem (\ref{3.2})-(\ref{3.6}) makes
no use, in fact, of the geometry $\partial\Omega_g$, it is defined in $D$ and
has admissible controls $g,u$ as defined in this section. We have:

\begin{proposition}\label{prop:3.1}
  The shape optimization problem (\ref{2.1}), (\ref{2.2}), (\ref{2.15}),
  (with $h\equiv 1$ in $D$) is equivalent with the constrained optimal control
  problem (\ref{3.2})-(\ref{3.6}).
\end{proposition}

\textit{Proof.} The cost computed via (\ref{2.15}) or (\ref{3.4}) are identical
since the solution of (\ref{2.9})-(\ref{2.11}) is a global parametrization
of $\partial \Omega_g$, on $[0,T_g]$. Moreover, it is clear that
$y_g|_{\Omega_g}$ satisfies (\ref{2.1}) and (\ref{3.5}), (\ref{3.6}) show
the $y_g$ also satisfies (\ref{2.2}). That is any admissible control
$g\in \mathcal{F}$ gives an admissible domain $\Omega_g$ and the two costs are
equal.

Conversely, starting from the admissible domain of type $\Omega=\Omega_g$,
$g\in \mathcal{F}$, one can extend the state system to $D$ by using an appropriate measurable
control $u$ , such that $g_+^2 u \in L^2(D)$ and the corresponding cost
is not changed via this transformation. This ends the proof.
\quad $\Box$

\begin{remark}\label{rem:3.2}
The Dirichlet boundary conditions (\ref{3.2}) on $\partial D$ may be replaced by any other
convenient boundary condition (for instance, simply supported plates in D). Here, we work
in $H^4(D) \cap H^2_0(D)$ that is a space of test functions, but in numerical experiments
other choices may be used. In \cite{MT2022} similar results are proved for second order state equations
and general cost functionals.
One can also work in $L^p(D)$, $p\geq 1$, according to \cite{Agmon1959}.
The constraints (\ref{3.5}), (\ref{3.6}) may be written as one, by addition.
\end{remark}

A standard technique in constrained optimization problems is the penalization
of the constraints in the cost functional, \cite{Bertsekas1999} ($\epsilon>0$
is ``small''):
\begin{eqnarray}
&&\min_{g,u}
\int_0^{T_g}  
\left[
  j\left(z_1(t),z_2(t),\Delta y_g(z_1(t),z_2(t))\right)
  +\frac{1}{\epsilon}\left| y_g\left(z_1(t),z_2(t)\right) \right|^2
\right.
\nonumber\\
&&
+\left.
  \frac{1}{\epsilon}
  \left| \nabla y_g\left(z_1(t),z_2(t)\right)\cdot
  \frac{\nabla g\left(z_1(t),z_2(t)\right)}{|\nabla g\left(z_1(t),z_2(t)\right)|}
  \right|^2
\right]
  |(z_1^\prime(t),z_2^\prime(t)) | dt
\label{3.7}
\end{eqnarray}
subject to (\ref{3.2})-(\ref{3.3}) and (\ref{2.9})-(\ref{2.11}). As explained in the end of
the previous section, the penalized  cost functional (\ref{3.7}) may be in fact a sum
if $\partial \Omega_g$ has more components, due to Thm. \ref{theo:2.1}. All the
considerations below may be easily adapted to such a situation. 

Both optimization problems (\ref{2.13})
and (\ref{3.4}) (with their state systems and restriction conditions if necessary)
may have no global optimal solutions due to the lack of compactness assumptions. We use instead just minimizing
sequences.

\begin{proposition}\label{prop:3.3}
  Let $j(\cdot, \cdot)$ be a positive Carath\'eodory function on $D\times\mathbb{R}$.
  Assume the $\mathcal{F}\subset \mathcal{C}^4(\overline{D})$ satisfies
  (\ref{2.6}), (\ref{2.7}) and $[y_n^\epsilon,g_n^\epsilon,u_n^\epsilon]$ is a
  minimizing sequence for the penalized problem
  (\ref{3.7}), (\ref{3.2}), (\ref{3.3}), (\ref{2.9})-(\ref{2.11}). Then, on a
  subsequence $n(m)$, the cost (\ref{3.4}) corresponding to
  $[\Omega_{g_{n(m)}^\epsilon}, y_{n(m)}^\epsilon]$ approaches some value less that
  $\inf \{  (\ref{3.2})-(\ref{3.6}), (\ref{2.9})-(\ref{2.11})\}$,
  (\ref{2.1}) is satisfied by $y_{n(m)}^\epsilon$ in $\Omega_{g_{n(m)}^\epsilon}$ and
  (\ref{2.2}) is satisfied with an error of order $\epsilon^{1/2}$.
\end{proposition}

\textit{Proof.}
The argument is similar to \cite{Tiba2018a}, \cite{MT2022}.  We take
$[y_{g_m}, \Omega_{g_m}]$ to be a minimizing sequence for the shape
optimization problem (\ref{2.1}), (\ref{2.2}), (\ref{2.15}) and we have
$y_{g_m}\in H^4(\Omega_{g_m})$ due to $\mathcal{F}\subset \mathcal{C}^4(\overline{D})$
and \cite{Agmon1959}. Then, we can construct the extension
$\tilde{y}_{g_m}\in H^4(D\setminus\overline{\Omega}_{g_m})
\cap H^2_0(D\setminus\overline{\Omega}_{g_m})$
and the control
\begin{equation}\label{3.8}
  u_{g_m}=\frac{\Delta\Delta \tilde{y}_{g_m} -f}{(g_m)_+^2}
  \quad\hbox{in }D\setminus\overline{\Omega}_{g_m},
\end{equation}
such that $u_{g_m}(g_m)_+^2\in L^2(D)$. This is admissible for the problem
(\ref{3.2})-(\ref{3.6}) and the corresponding state is obtained by the concatenation
of $y_{g_m}$, $\tilde{y}_{g_m}$. The construction is also valid for $\Omega_{g_m}$
not necessarily simply connected. The cost (\ref{2.15}) majorizes the one
in (\ref{3.4}) or (\ref{3.7}) (the penalization term is null) and
we get the inequality
\begin{eqnarray}
&&
\int_0^{T_g}  
\left[
  j\left(\mathbf{z}_{n(m)}^\epsilon,\Delta y_{n(m)}^\epsilon(\mathbf{z}_{n(m)}^\epsilon)\right)
  +\frac{1}{\epsilon}\left| y_{n(m)}^\epsilon\left(\mathbf{z}_{n(m)}^\epsilon\right) \right|^2
\right.
\nonumber\\
&&
+\left.
  \frac{1}{\epsilon}
  \left| \nabla y_{n(m)}^\epsilon\left(\mathbf{z}_{n(m)}^\epsilon\right)\cdot
  \frac{\nabla g\left(\mathbf{z}_{n(m)}^\epsilon\right)}
       {\left|\nabla g\left(\mathbf{z}_{n(m)}^\epsilon\right)\right|}
  \right|^2
\right]
| (\mathbf{z}_{n(m)}^\epsilon)^\prime | dt
\nonumber\\
&\leq &
\int_{\partial \Omega_{g_m}} j(\sigma,\Delta y_m(\sigma))d\sigma
\rightarrow
\inf\{ (\ref{2.1}), (\ref{2.2}), (\ref{2.15}\}.
\label{3.9}
\end{eqnarray}
In (\ref{3.9}), the index $n(m)$ is big enough in order to have the inequality
satisfied (due to the admissibility of control pair $[g_m,u_m]$) and
$\mathbf{z}_{n(m)}^\epsilon$ is the solution of (\ref{2.9})-(\ref{2.11}) associated to
$g_{n(m)}^\epsilon$.

Inequality (\ref{3.9}), due to the positivity of $j(\cdot,\cdot)$ gives the last
statement of the conclusions. The equation (\ref{2.1}) is satisfied by
$y_{n(m)}^\epsilon$ by the construction of (\ref{3.3}) and the minimizing property
is again a consequence of (\ref{3.4}).

\begin{remark}\label{rem:3.4}
  In this section we use the above extension of the state system together with the classical
  constraints penalization technique in the cost functional,
  \cite{Bertsekas1999}, that is a very general approach. Another variant is the penalization
   of the state system (\ref{2.1}), (\ref{2.2}). This was
  introduced in  free boundary problems by 
  Kawarada and  Natori \cite{Kawarada1981} and used in shape optimization in
  \cite{NP_Tiba2009}, \cite{N_Tiba2012}. It was also applied to second
  order elliptic equations, with Dirichlet, respectively Neumann boundary conditions
  by \cite{Zhou2014}, \cite{Zhou2017}. It is to be noticed that the penalization technique has a
  rather different structure in the two cases. The problem that we study
  in this paper is a Dirichlet problem for the biharmonic operator
  and one may say that it ``includes'' in the boundary conditions both conditions
  (Dirichlet and Neumann) from second order problems. Therefore, it seems
  unclear, in this case, how to apply the penalization method for the state
  system and the only possible approach seems to be via the Hamiltonian
  system (\ref{2.9})-(\ref{2.11}) and the penalization of the cost functional, that we
	use here.
\end{remark}

We consider now functionals $g+\lambda r$, $u+\lambda v$ with $\lambda\in \mathbb{R}$,
$g,r\in\mathcal{F}\subset\mathcal{C}^2(\overline{D})$ and $u,v\in L^2(D)$.
Notice that this ensures, for instance, $(g+\lambda r)_+^2(u+\lambda v)\in L^2(D)$,
as required. Our aim is to compute the directional derivative of the cost functional
(\ref{3.7}) to be used in the numerical applications from the next section.
We have

\begin{proposition}\label{prop:3.5}
The limit $\mathbf{w}=[w_1,w_2]=\lim_{\lambda\rightarrow 0}
\frac{\mathbf{z}_{\lambda} - \mathbf{z}}{\lambda}$
exists in $\mathcal{C}^2([0,T_g])^2$ and satisfies
\begin{eqnarray}
w_1^\prime & = & -\nabla\partial_2 g(\mathbf{z})\cdot \mathbf{w}
-\partial_2 r(\mathbf{z}),\quad\hbox{in } [0,T_g],
\label{3.10}\\
w_2^\prime & = & \nabla\partial_1 g(\mathbf{z})\cdot \mathbf{w}
+\partial_1 r(\mathbf{z}),\quad\hbox{in } [0,T_g],
\label{3.11}\\
\mathbf{w}(0) & = & [w_1(0),w_2(0)] = [0,0],
\label{3.12}
\end{eqnarray}
where $\mathbf{z}$ is the solution of (\ref{2.9})-(\ref{2.11}) and
$\mathbf{z}_\lambda$ is the solution of the perturbed variant of (\ref{2.9})-(\ref{2.11})
with $g+\lambda r$ instead of $g$.
\end{proposition}

The proof of (\ref{3.10})-(\ref{3.11}) is standard and can be found in 
\cite{Tiba2013}. Notice that $g+\lambda r$ satisfies (\ref{2.6}), (\ref{2.7}) for
$| \lambda |$ small, by the application of the Weierstrass theorem and, consequently,
$\mathbf{z}_\lambda$ is periodic by Thm. \ref{theo:2.1} and the ratio
$\frac{\mathbf{z}_{\lambda} - \mathbf{z}}{\lambda}$ is well defined on $[0,T_g]$.

\begin{proposition}\label{prop:3.6}
Denote by $q=\lim_{\lambda\rightarrow 0}\frac{y_{\lambda} - y_g}{\lambda}$ that exists
in $H^4(D)$ strong, where $y_{\lambda}=y_{g+\lambda r}$ is the solution of the perturbed
equation (\ref{3.3}) with $(g+\lambda r)_+^2(u+\lambda v)$ in the right-hand side.
Then $q\in H^4(D)\cap H_0^2(D)$ satisfies:
\begin{eqnarray}
\Delta\Delta q & = & g_+^2v+2g_+ur \hbox{ in }D,
\label{3.13}\\
q & = & \frac{\partial q}{\partial \mathbf{n}}=0\hbox{ on }\partial D.
\label{3.14}
\end{eqnarray}
\end{proposition}

\textit{Proof.}
We indicate a sketch of the argument since the steps are quite standard. For second order
equations, we quote \cite{MT2019}. Subtracting (\ref{3.3}) and its perturbed variant
we get,
\begin{eqnarray}
  \Delta\Delta \frac{y_{\lambda} - y_g}{\lambda} & = & (g+\lambda r)_+^2 v
  +\frac{1}{\lambda}\left[ (g+\lambda r)_+^2 -g_+^2\right]u
  \hbox{ in }D,
\label{3.15}\\
\frac{y_{\lambda} - y_g}{\lambda} & = &
\frac{\partial}{\partial \mathbf{n}} \left(\frac{y_{\lambda} - y_g}{\lambda}\right)
=0\hbox{ on }\partial D.
\label{3.16}
\end{eqnarray}

Obviously, the right-hand side in (\ref{3.15}) has the limit appearing in (\ref{3.13}),
in the space $L^2(D)$. The estimate for strong solutions of Dirichlet problems
for biharmonic equations \cite{Agmon1959}, \S 8, show that
$\left\{ \frac{y_{\lambda} - y_g}{\lambda}\right\}$ is bounded in  $H^4(D)\cap H_0^2(D)$.
Then, the limit $q \in H^4(D)\cap H_0^2(D)$ exists on a subsequence, in the weak
topology of $H^4(D)$, and satisfies (\ref{3.13}), (\ref{3.14}). Since the convergence in the 
right-hand side of (\ref{3.15}) is strong, the same is valid for
$\left\{ \frac{y_{\lambda} - y_g}{\lambda}\right\}$ in $H^4(D)$. The uniqueness
of the limit $q$ shows that the convergence is valid without taking subsequence.
\quad$\Box$

To study the differentiability properties of the penalized cost function (\ref{3.7}), we also assume
$f \in H^1(D)$. By properties of the positive part, we get that $g^2_+ \in
W^{1,\infty}(D)$ and $g^2_+ u \in H^1(D)$ if $u \in H^1(D)$ and the solution of (\ref{3.2}), (\ref{3.3}) satisfies, by Ch.9, \cite{Agmon1959}, that $y_g \in H^5(D) \subset C^3(\overline{D})$ (due to the Sobolev
theorem).
\begin{theorem}\label{theo:3.7}
  Under the above conditions, assume that $j(\cdot,\cdot)\in\mathcal{C}^1(\mathbb{R}^3)$.
  Then, the directional derivative of the approximating performance
  index (\ref{3.7}), in the point $[g,u]\in \mathcal{F}\times L^2(D)$ and in the
  direction $[r,v]\in \mathcal{F}\times L^2(D)$, is given by:
\begin{eqnarray}
&&
\theta(g,r)\left[
j(\mathbf{x}_0,\Delta y_g(\mathbf{x}_0))
+\frac{1}{\epsilon} \left| y_g(\mathbf{x}_0) \right|^2
+\frac{1}{\epsilon} \left|
\nabla y_g(\mathbf{x}_0) \cdot
\frac{\nabla g(\mathbf{x}_0)}{|\nabla g(\mathbf{x}_0)|}
\right|^2
\right] | \nabla g(\mathbf{x}_0)|
\nonumber\\
&+& 
\int_0^{T_g}\left[
\nabla_1 j\left(\mathbf{z}(t),\Delta y_g(\mathbf{z}(t))\right)
\cdot \mathbf{w}(t)
\right] |\mathbf{z}^\prime(t)| dt
\nonumber\\
&+&
\int_0^{T_g}
\partial_2 j\left(\mathbf{z}(t),\Delta y_g(\mathbf{z}(t))\right)
\left[
\nabla \Delta y_g(\mathbf{z}(t))\cdot \mathbf{w}(t)
+\Delta q(\mathbf{z}(t))
\right]
|\mathbf{z}^\prime(t)| dt
\nonumber\\
&+&
\int_0^{T_g}
j\left(\mathbf{z}(t),\Delta y_g(\mathbf{z}(t))\right)
\frac{\mathbf{z}^\prime(t)\cdot \mathbf{w}^\prime(t)}
{|\mathbf{z}^\prime(t)|}
dt
\nonumber\\
&+&
\frac{2}{\epsilon}\int_0^{T_g}
\left[
  y_g(\mathbf{z}(t))\nabla y_g(\mathbf{z}(t))\cdot \mathbf{w}(t)
  +y_g(\mathbf{z}(t))q(\mathbf{z}(t))
\right]
|\mathbf{z}^\prime(t)| dt
\nonumber\\
&+&
\frac{1}{\epsilon}\int_0^{T_g}
\left| y_g(\mathbf{z}(t))\right|^2
\frac{\mathbf{z}^\prime(t)\cdot \mathbf{w}^\prime(t)}
{|\mathbf{z}^\prime(t)|}
dt
\nonumber\\
&+&\frac{2}{\epsilon}
\int_0^{T_g}
\left( \nabla y_g(\mathbf{z}(t))\cdot
\frac{\nabla g(\mathbf{z}(t))}{|\nabla g(\mathbf{z}(t))|}
\right)
\frac{\nabla r(\mathbf{z}(t))}{|\nabla g(\mathbf{z}(t))|}
\cdot\nabla y_g(\mathbf{z}(t))
|\mathbf{z}^\prime(t)|
dt
\nonumber\\
&+&
\frac{2}{\epsilon}
\int_0^{T_g}
\left(
\nabla y_g(\mathbf{z}(t))\cdot
\frac{\nabla g(\mathbf{z}(t))}{|\nabla g(\mathbf{z}(t))|}
\right)
\left[
\left(H\,y_g(\mathbf{z}(t))\right)\mathbf{w}(t)
\right]
\cdot
\frac{\nabla g(\mathbf{z}(t))}{|\nabla g(\mathbf{z}(t))|}
|\mathbf{z}^\prime(t)|
dt
\nonumber\\
&+&\frac{2}{\epsilon}
\int_0^{T_g}
\left(
\nabla y_g(\mathbf{z}(t))\cdot
\frac{\nabla g(\mathbf{z}(t))}{|\nabla g(\mathbf{z}(t))|}
\right)
\nabla q(\mathbf{z}(t))
\cdot
\frac{\nabla g(\mathbf{z}(t))}{|\nabla g(\mathbf{z}(t))|}
|\mathbf{z}^\prime(t)|
dt
\nonumber\\
&+&\frac{2}{\epsilon}
\int_0^{T_g}
\left(
\nabla y_g(\mathbf{z}(t))\cdot
\frac{\nabla g(\mathbf{z}(t))}{|\nabla g(\mathbf{z}(t))|}
\right)
\nabla y_g(\mathbf{z}(t))\cdot
\left[
\frac{\left(H\,g(\mathbf{z}(t))\right)\mathbf{w}(t)}
{|\nabla g(\mathbf{z}(t))|}
\right.  
\nonumber\\
&&-\left.
\frac{\nabla g(\mathbf{z}(t))}{|\nabla g(\mathbf{z}(t))|^3}
\Big(
\nabla g(\mathbf{z}(t)) \cdot \nabla r(\mathbf{z}(t))
+\nabla g(\mathbf{z}(t)) \cdot
\left(H\,g(\mathbf{z}(t))\right)\mathbf{w}(t)
\Big)
\right]|\mathbf{z}^\prime(t)|
dt
\nonumber\\
&+&
\frac{1}{\epsilon}
\int_0^{T_g}
\left[
\nabla y_g(\mathbf{z}(t))\cdot
\frac{\nabla g(\mathbf{z}(t))}{|\nabla g(\mathbf{z}(t))|}
\right]^2
\frac{\mathbf{z}^\prime(t)\cdot \mathbf{w}^\prime(t)}
{|\mathbf{z}^\prime(t)|}
dt,
\label{3.17}
\end{eqnarray}
where $\nabla_1 j$ is the gradient of $j(\cdot,\cdot)$
with respect to the vector $\mathbf{z}$ and 
$\partial_2 j$ is the partial derivative with respect to the last variable,
$H\,y_g$ is the Hessian matrix of $y_g\in H^4(D)$ and $q$, $w$ are defined
in Propositions \ref{prop:3.5}, \ref{prop:3.6}, etc.
\end{theorem}

\textbf{Proof.} We indicate just some computations due the similarity with
\cite{MT2022} \S 4. We discuss first the term:
\begin{eqnarray*}
&&\frac{1}{\lambda}
\int_{T_g}^{T_\lambda}
\left[ 
  j\left(\mathbf{z}_\lambda(t),\Delta y_\lambda(\mathbf{z}_\lambda(t))\right)
  +\frac{1}{\epsilon} \left| y_\lambda(\mathbf{z}_\lambda(t)) \right|^2
  \right.
\nonumber\\
&&  
\left.  
  +\frac{1}{\epsilon}
  \left|
    \nabla y_\lambda(\mathbf{z}_\lambda(t)) \cdot
    \frac{\nabla (g+\lambda r)(\mathbf{z}_\lambda(t))}{|\nabla (g+\lambda r)(\mathbf{z}_\lambda(t))|}
  \right|^2
\right]
|\mathbf{z}_\lambda^\prime(t)|
dt
\nonumber\\
& = &
\frac{T_\lambda-T_g}{\lambda}
\left[
j\left(\mathbf{z}_\lambda(\tau_\lambda),
\Delta y_\lambda(\mathbf{z}_\lambda(\tau_\lambda))\right)
+\frac{1}{\epsilon} \left| y_\lambda(\mathbf{z}_\lambda(\tau_\lambda)) \right|^2
  \right.
\nonumber\\
&&  
\left.
  +\frac{1}{\epsilon}
  \left|
    \nabla y_\lambda(\mathbf{z}_\lambda(\tau_\lambda)) \cdot
    \frac{\nabla (g+\lambda r)(\mathbf{z}_\lambda(\tau_\lambda))}
         {|\nabla (g+\lambda r)(\mathbf{z}_\lambda(\tau_\lambda))|}
  \right|^2
\right]
|\mathbf{z}_\lambda^\prime(\tau_\lambda)|
\nonumber\\
& \rightarrow &
\theta(g,r)
\left[
  j(\mathbf{x}_0,\Delta y_g(\mathbf{x}_0))
  +\frac{1}{\epsilon} \left| y_g(\mathbf{x}_0) \right|^2
  +\frac{1}{\epsilon} \left|
\nabla y_g(\mathbf{x}_0) \cdot
\frac{\nabla g(\mathbf{x}_0)}{|\nabla g(\mathbf{x}_0)|}
\right|^2
\right]  
|\mathbf{z}^\prime(T_g)|,
\end{eqnarray*}
where $\tau_\lambda$ is some intermediary point between $T_g$ and $T_\lambda$ and
we use the periodicity of $\mathbf{z}_\lambda(\cdot)$, $\mathbf{z}(\cdot)$ and the system
(\ref{2.9})-(\ref{2.11}) for the values of $\mathbf{z}_\lambda^\prime(\cdot)$.
We also use the continuity of $j$, $\mathbf{z}_\lambda$, $y_\lambda$, $\nabla y_\lambda$,
$\nabla g$, $\Delta y_\lambda$, their convergence properties and Proposition \ref{prop:2.2}.
Notice that $y_g(\mathbf{x}_0)$, $\nabla y_g(\mathbf{x}_0)\cdot
\frac{\nabla g(\mathbf{x}_0)}{|\nabla g(\mathbf{x}_0)|}
=\frac{\partial g(\mathbf{x}_0)}{\partial \mathbf{n}}$ are not necessarily null since the
constraints (\ref{3.5}), (\ref{3.6}) may not be satisfied in the penalized problem.

We continue with the terms:
\begin{eqnarray*}
&&
\frac{1}{\lambda}
\int_0^{T_g}
\bigg\{
\left[
j\left(\mathbf{z}_\lambda(t),\Delta y_\lambda(\mathbf{z}_\lambda(t))\right)
+\frac{1}{\epsilon}
\left| y_\lambda(\mathbf{z}_\lambda(t)) \right|^2
\right.
\nonumber\\
&&
\left.
+\frac{1}{\epsilon}
\left|
\nabla y_\lambda(\mathbf{z}_\lambda(t))\cdot
\frac{\nabla (g+\lambda r)(\mathbf{z}_\lambda(t))}
{|\nabla (g+\lambda r)(\mathbf{z}_\lambda(t))|} 
\right|^2
\right]
|\mathbf{z}_\lambda^\prime(t)|
\nonumber\\
&&
-
\left[
j\left(\mathbf{z}(t),\Delta y_g(\mathbf{z}(t))\right)
+\frac{1}{\epsilon} \left| y_g(\mathbf{z}(t)) \right|^2
\right.
\nonumber\\
&&
\left.
+\frac{1}{\epsilon}
\left|
\nabla y_g(\mathbf{z}(t))\cdot
\frac{\nabla g(\mathbf{z}(t))}
{|\nabla g(\mathbf{z}(t))|} 
\right|^2
\right]
|\mathbf{z}^\prime(t)|
\bigg\}
dt.
\end{eqnarray*}

For the part corresponding to $j(\cdot,\cdot)$ and the first penalized term, we get
\begin{eqnarray*}
&&\int_0^{T_g}
\left[ \nabla_1 j\left(\mathbf{z}(t),\Delta y_g(\mathbf{z}(t))\right)
\cdot \mathbf{w}(t) 
+\partial_2 j\left(\mathbf{z}(t),\Delta y_g(\mathbf{z}(t))\right)
\nabla \Delta y_g(\mathbf{z}(t)) \cdot \mathbf{w}(t)
\right]
|\mathbf{z}^\prime(t)|
dt
\nonumber\\
&&+
\int_0^{T_g}
\partial_2 j\left(\mathbf{z}(t),\Delta y_g(\mathbf{z}(t))\right)
\Delta q(\mathbf{z}(t))
|\mathbf{z}^\prime(t)|
dt
+\int_0^{T_g}
j\left(\mathbf{z}(t),\Delta y_g(\mathbf{z}(t))\right)
\frac{\mathbf{z}^\prime(t)\cdot \mathbf{w}^\prime(t)}
     {|\mathbf{z}^\prime(t)|}
dt
\nonumber\\
&&+
\frac{2}{\epsilon}\int_0^{T_g} y_g(\mathbf{z}(t))
\nabla  y_g(\mathbf{z}(t))\cdot \mathbf{w}(t)|\mathbf{z}^\prime(t)|
dt
+
\frac{2}{\epsilon}\int_0^{T_g} y_g(\mathbf{z}(t))
q(\mathbf{z}(t))|\mathbf{z}^\prime(t)|
dt
\nonumber\\
&&+\frac{1}{\epsilon}\int_0^{T_g}
\left| y_g(\mathbf{z}(t)) \right|^2
\frac{\mathbf{z}^\prime(t)\cdot \mathbf{w}^\prime(t)}
     {|\mathbf{z}^\prime(t)|}
dt
\end{eqnarray*}
and the arguments are very similar with \cite{MT2022}, \S 4.
The handling of the penalization terms containing
$\nabla y_\lambda$, $\nabla y_g$ is performed as in \cite{MT2022}, for $\alpha\equiv 0$.
\quad$\Box$

\begin{remark}\label{rem:3.8}
  It is possible to consider more general cost functionals, as in \cite{MT2022}.
  We have studied just the case of boundary observation given in
  (\ref{2.15}) since we consider it the most significant and difficult.
\end{remark}

\section{Finite element approximation and descent directions algorithm}
\setcounter{equation}{0}

We assume that $D$ is a polyhedral domain.
We denote by $\mathcal{T}_h$ a triangulation of $D$, with $h$ the mesh size.
The finite element approximation of $g$, $u$, and $y_g$ are denoted
$g_h$, $u_h$, $y_h$, respectively. We have $g_h\in \mathbb{W}_h^g$,
$u_h\in \mathbb{W}_h^u$ and $y_g\in \mathbb{V}_h$, where
\begin{eqnarray*}
\mathbb{W}_h^g & = & \{ g_h\in \mathcal{C}(\overline{D});
\ {g_h}_{|T} \in \mathbb{P}_3(T),\ \forall T \in \mathcal{T}_h  \}\\
\mathbb{W}_h^u & = & \{ u_h\in \mathcal{C}(\overline{D});
\ {u_h}_{|T} \in \mathbb{P}_1(T),\ \forall T \in \mathcal{T}_h  \}
\end{eqnarray*}
of dimensions $n_g=card(\mathbb{W}_h^g)$, $n_u=card(\mathbb{W}_h^u)$
and basis $\{ \phi_i^g;\ 1\leq i\leq n_g\}$,
$\{ \phi_i^u;\ 1\leq i\leq n_u\}$, respectively.
We precise that $g_h$ is not necessarily in $\mathcal{C}^1(\overline{D})$, but
in each triangle $T$ of $\mathcal{T}_h$, the first and second order derivatives
of ${g_h}_{|T}$ are well defined.

The vectorial space $\mathbb{V}_h$ is based on the Hsieh, Clough and Tocher (HCT)
finite element, where each triangle $T$ is subdivided in 3 triangles $T_i$,
the vertices of $T_i$ are the barycenter of $T$ and two vertices of $T$,
$y_h{}_{|T_i} \in \mathbb{P}_3(T_i)$, $i=1,2,3$, the value of the function and its first
derivatives are continuous at each vertex of $T$, the normal derivative
is continuous at the mid points of the 3 sides of $T$,
12 degrees of freedom, globally of class $\mathcal{C}^1$,
see \cite{Dautray1990}, Ch. XII, section 4.
More precisely, we set
$$
\mathbb{V}_h  =  \{ y_h\in \mathcal{C}^1(\overline{D});
\ {y_h}_{|T} \hbox{ HCT finite element},\ \forall T \in \mathcal{T}_h,
\ y_h=\frac{\partial y_h}{\partial \mathbf{n}}=0\hbox{ on }\partial D\}.
$$

We can write
$$
g_h(\mathbf{x})=\sum_{i=1}^{n_g} G_i \phi_i^g(\mathbf{x}),\quad
u_h(\mathbf{x})=\sum_{i=1}^{n_u} U_i \phi_i^u(\mathbf{x})
$$
where $G=(G_i)_{1\leq i\leq n_g} \in \mathbb{R}^{n_g}$ and
$U=(U_i)_{1\leq i\leq n_u} \in \mathbb{R}^{n_u}$.
The discrete version of (\ref{3.7}) is
\begin{eqnarray}
  \min_{G,U} \mathcal{J}(G,U)
&=&
\int_{\partial\Omega_h}
j\left(\sigma,\Delta y_h(\sigma)\right)
d\sigma
+\frac{1}{\epsilon} \int_{\partial\Omega_h}
\left| y_h\left(\sigma\right) \right|^2
d\sigma
\nonumber\\
&&
+ \frac{1}{\epsilon}\int_{\partial\Omega_h}
  \left| \nabla y_h\left(\sigma\right)\cdot
  \frac{\nabla g_h\left(\sigma\right)}{|\nabla g_h\left(\sigma\right)|}
  \right|^2
d\sigma 
\label{4.1}
\end{eqnarray}
where  $y_h\in \mathbb{V}_h$ is the solution of
\begin{equation}\label{4.2}
\int_D \Delta y_h  \Delta \varphi_h \, d\mathbf{x}
= \int_D \left( f\varphi_h + (g_h)_+^2 u_h \varphi_h \right) d\mathbf{x},
\qquad \forall \varphi_h \in \mathbb{V}_h
\end{equation}
and $\partial\Omega_h$ is composed by closed polygonal lines
approximating $\{\mathbf{x}\in D;\ g_h(\mathbf{x})=0\}$.
More precisely, the connected component indexed ``c'' of $\partial\Omega_h$ is the union
on $n_c$ segments, i.e. $\cup_{\ell=1}^{n_c} [S_{\ell-1}^c,S_{\ell}^c]$ and for all
$\ell$, there exists $T\in \mathcal{T}_h$, such that
$S_{\ell-1}^c,S_{\ell}^c\in \partial T$. In other words, the vertices $S_{\ell}^c$ of the
polygonal lines are on the sides of the triangulation $\mathcal{T}_h$.
Then, $\nabla g_h$ and $\nabla y_h$ are well defined on the segment
$[S_{\ell-1}^c,S_{\ell}^c]\subset T$.
For $y_h\in \mathbb{V}_h$, we have only ${y_h}_{|T}\in \mathcal{C}^1(\overline{T})$, so
we have to explain the meaning of $\Delta y_h$ appearing in
(\ref{4.1}), (\ref{4.2}). In the finite element software FreeFem++ \cite{freefem++},
there exist the discrete second order derivative operators \texttt{dxx}, \texttt{dyy}
and $\Delta y_h:=dxx(y_h)+dyy(y_h)$ is a $\mathbb{P}_1$ finite element function
continuous in $D$ approaching $\Delta y_g$. We note that these operators use interpolation, too.
Finally, in (\ref{4.1}), the discrete objective function is the sum of integrals
of continuous functions over segments.

We denote by $t_1$ the first term of (\ref{4.1}),
$$
t_2=\int_{\partial\Omega_h}
\left| y_h\left(\sigma\right) \right|^2
d\sigma, \quad
t_3=\int_{\partial\Omega_h}
  \left| \nabla y_h\left(\sigma\right)\cdot
  \frac{\nabla g_h\left(\sigma\right)}{|\nabla g_h\left(\sigma\right)|}
  \right|^2
d\sigma 
$$
and we get
$\mathcal{J}(G,U)=t_1 +\frac{1}{\epsilon}t_2 +\frac{1}{\epsilon}t_3$.

\medskip
\textbf{Algorithm.}

We employ a descent direction method
$$
(G^{k+1},U^{k+1})=(G^k,U^k)+\lambda_k (R^k,V^k),
$$
where $\lambda_k >0$ is obtained by line search
$$
\lambda_k \in \arg\min_{\lambda >0}
\mathcal{J}\left((G^k,U^k)+\lambda (R^k,V^k)\right)
$$
and $(R^k,V^k)$ is a descent direction, i.e. $d\mathcal{J}_{(G^k,U^k)}(R^k,V^k) <0$,
where $d\mathcal{J}_{(G,U)}(R,V)$ is a discrete version of (\ref{3.17}), detailed
in the following and given by the formula (\ref{4.15}).
We stop the algorithm when
$| \mathcal{J}(G^{k+1},U^{k+1}) - \mathcal{J}(G^k,U^k)| < tol$.

\medskip

Now, we discuss the approximation of the directional derivative given by (\ref{3.17}).
We can replace $\mathbf{w}^\prime$ from (\ref{3.17}) using (\ref{3.10}), (\ref{3.11}),
then the directional derivative is equal to $\Gamma_q+\Gamma_r+\Gamma_w$, where
$\Gamma_q$ is the sum of the three terms containing $q$,
$\Gamma_r$ is the sum of the five terms containing $r$ and $\Gamma_w$ the remaining
terms containing $\mathbf{w}$, including the first term of (\ref{3.17}).

\medskip
\textbf{Approximation of the terms containing $q$}

We have
\begin{eqnarray*}
\Gamma_q&=&\int_0^{T_g}
\partial_2 j\left(\mathbf{z}(t),\Delta y_g(\mathbf{z}(t))\right)\Delta q(\mathbf{z}(t))
|\mathbf{z}^\prime(t)| dt
+\frac{2}{\epsilon} \int_0^{T_g}
y_g\left(\mathbf{z}(t)\right)q(\mathbf{z}(t))
|\mathbf{z}^\prime(t)|
dt
\nonumber\\
&&
+ \frac{2}{\epsilon}\int_0^{T_g}
\left(
\nabla y_g\left(\mathbf{z}(t)\right)\cdot
\frac{\nabla g\left(\mathbf{z}(t)\right)}{|\nabla g\left(\mathbf{z}(t)\right)|}
\right)
\nabla q(\mathbf{z}(t))\cdot 
\frac{\nabla g\left(\mathbf{z}(t)\right)}{|\nabla g\left(\mathbf{z}(t)\right)|}
|\mathbf{z}^\prime(t)|
dt 
\end{eqnarray*}
and we can write it as an integral over $\partial\Omega_g$,
\begin{eqnarray*}
\Gamma_q(y_g,g,q)&=&\int_{\partial\Omega_g}
\partial_2 j\left(\sigma,\Delta y_g(\sigma)\right)\Delta q(\sigma)
d\sigma
+\frac{2}{\epsilon} \int_{\partial\Omega_g}
y_g\left(\sigma\right)q(\sigma)
d\sigma
\nonumber\\
&&
+ \frac{2}{\epsilon}\int_{\partial\Omega_g}
\left(
\nabla y_g\left(\sigma\right)\cdot
\frac{\nabla g\left(\sigma\right)}{|\nabla g\left(\sigma\right)|}
\right)
\nabla q(\sigma)\cdot 
\frac{\nabla g\left(\sigma\right)}{|\nabla g\left(\sigma\right)|}
d\sigma.
\end{eqnarray*}

Let us introduce the discrete weak formulation of (\ref{3.13})-(\ref{3.14}): 
for $r_h\in \mathbb{V}_h^g$ and $v_h\in \mathbb{V}_h^u$,
find $q_h\in \mathbb{V}_h$ such that
\begin{equation}\label{4.3}
\int_D \Delta q_h  \Delta \varphi_h \, d\mathbf{x}
= \int_D  \left( (g_h)_+^2 v_h \varphi_h + 2(g_h)_+ u_h r_h \varphi_h \right) d\mathbf{x},
\qquad \forall \varphi_h \in \mathbb{V}_h
\end{equation}
and an auxiliary problem:
find $p_h\in \mathbb{V}_h$ such that
\begin{eqnarray}
\int_D \Delta \varphi_h  \Delta p_h \, d\mathbf{x}
&= &\int_{\partial\Omega_h}
\partial_2 j\left(\sigma,\Delta y_h(\sigma)\right)\Delta\varphi_h(\sigma)
d\sigma
+\frac{2}{\epsilon} \int_{\partial\Omega_h}
y_h\left(\sigma\right)\varphi_h(\sigma)
d\sigma
\nonumber\\
&&
+ \frac{2}{\epsilon}\int_{\partial\Omega_h}
\left(
\nabla y_h\left(\sigma\right)\cdot
\frac{\nabla g_h\left(\sigma\right)}{|\nabla g_h\left(\sigma\right)|}
\right)
\nabla \varphi_h(\sigma)\cdot 
\frac{\nabla g_h\left(\sigma\right)}{|\nabla g_h\left(\sigma\right)|}
d\sigma, 
\label{4.4}
\end{eqnarray}
for all $\varphi_h \in \mathbb{V}_h$.

Putting $\varphi_h=p_h$ in (\ref{4.3}) and $\varphi_h=q_h$ in (\ref{4.4}), we get
\begin{eqnarray}
&&
\int_D \left( (g_h)_+^2 p_h v_h  + 2(g_h)_+ u_h p_h r_h \right)
d\mathbf{x}
=
\int_{\partial\Omega_h}
\partial_2 j\left(\sigma,\Delta y_h(\sigma)\right)\Delta q_h(\sigma)
d\sigma
\nonumber\\
&&
+\frac{2}{\epsilon} \int_{\partial\Omega_h}
y_h\left(\sigma\right) q_h(\sigma)
d\sigma
\nonumber\\
&&
+\frac{2}{\epsilon}\int_{\partial\Omega_h}
\left(
\nabla y_h\left(\sigma\right)\cdot
\frac{\nabla g_h\left(\sigma\right)}{|\nabla g_h\left(\sigma\right)|}
\right)
\nabla q_h(\sigma)\cdot 
\frac{\nabla g_h\left(\sigma\right)}{|\nabla g_h\left(\sigma\right)|}
d\sigma
\label{4.5}
\end{eqnarray}
and the right-side hand of (\ref{4.5}) corresponds to $\Gamma_q(y_h,g_h,q_h)$
the approximation of $\Gamma_q$, the sum of
the terms of (\ref{3.17}) containing $q$.

Given $g_h\in\mathbb{W}_h^g$ and  $u_h\in\mathbb{W}_h^u$,  
let $y_h, p_h \in \mathbb{V}_h$ solutions of (\ref{4.3}), (\ref{4.4}),
respectively. We set
$$
r_h(\mathbf{x})=\sum_{i=1}^{n_g} R_i \phi_i^g(\mathbf{x}),\quad
R=(R_i)_{1\leq i\leq n_g},\quad
v_h(\mathbf{x})=\sum_{i=1}^{n_u} V_i \phi_i^u(\mathbf{x}),\quad
V=(V_i)_{1\leq i\leq n_u}
$$
and $R^q=(R_i^q)_{1\leq i\leq n_g} \in \mathbb{R}^{n_g}$,
$V^q=(V_i^q)_{1\leq i\leq n_u} \in \mathbb{R}^{n_u}$, given by
\begin{equation}\label{4.6}
R_i^q=-\int_D   2(g_h)_+ u_h p_h \phi_i^g 
\, d\mathbf{x},\ 1\leq i\leq n_g\quad
V_i^q=-\int_D  (g_h)_+^2 p_h \phi_i^u
\, d\mathbf{x}, 1\leq i\leq n_u.
\end{equation}
We obtain that
\begin{eqnarray*}
\Gamma_q(y_h,g_h,q_h)  
&=&\int_D \left( (g_h)_+^2 p_h v_h + 2(g_h)_+ u_h p_h r_h \right)
d\mathbf{x}\\
&=&\sum_{i=1}^{n_u} \int_D  (g_h)_+^2 p_h V_i \phi_i^u \, d\mathbf{x}
+\sum_{i=1}^{n_g} \int_D 2(g_h)_+ u_h p_h R_i \phi_i^g\, d\mathbf{x}\\
&=&\sum_{i=1}^{n_u} V_i \int_D  (g_h)_+^2 p_h \phi_i^u \, d\mathbf{x}
+\sum_{i=1}^{n_g} R_i \int_D 2(g_h)_+ u_h p_h  \phi_i^g\, d\mathbf{x}\\
&=&-\sum_{i=1}^{n_u} V_i^q V_i  -\sum_{i=1}^{n_g} R_i^q R_i
=-\langle V^q,V\rangle - \langle R^q,R\rangle.
\end{eqnarray*}
In particular, for $(R,V)=(R^q,V^q)$, we get that the approximation of $\Gamma_q$
is $-\|R^q\|^2 -\|V^q\|^2\leq 0$. If $p_h$ does not vanish in
$\{\mathbf{x}\in D;\  (g_h)_+(\mathbf{x}) >0\}$, then $V^q\neq 0$, and the above
inequality is strict.

\medskip
\textbf{Approximation of the terms containing $r$}

Using (\ref{2.9}), (\ref{2.10}), we get
\begin{eqnarray*}
\Gamma_r&=&
\frac{2}{\epsilon}
\int_0^{T_g}
\left( \nabla y_g(\mathbf{z}(t))\cdot
\frac{\nabla g(\mathbf{z}(t))}{|\nabla g(\mathbf{z}(t))|}
\right)
\frac{\nabla r(\mathbf{z}(t))}{|\nabla g(\mathbf{z}(t))|}
\cdot\nabla y_g(\mathbf{z}(t))
|\mathbf{z}^\prime(t)|  
dt
\nonumber\\
&-&
\frac{2}{\epsilon}
\int_0^{T_g}
\left(
\nabla y_g(\mathbf{z}(t))\cdot
\frac{\nabla g(\mathbf{z}(t))}{|\nabla g(\mathbf{z}(t))|}
\right)^2
\frac{1}{|\nabla g(\mathbf{z}(t))|^2}
\Big(
\nabla g(\mathbf{z}(t)) \cdot \nabla r(\mathbf{z}(t))
\Big)
|\mathbf{z}^\prime(t)|
dt
\nonumber\\
&+&
\int_0^{T_g}
j\left(\mathbf{z}(t),\Delta y_g(\mathbf{z}(t))\right)
\frac{1}{|\nabla g(\mathbf{z}(t))|}
\nabla g(\mathbf{z}(t)) \cdot \nabla r(\mathbf{z}(t))
dt
\nonumber\\
&+&
\frac{1}{\epsilon}
\int_0^{T_g}
\left[
\nabla y_g(\mathbf{z}(t))\cdot
\frac{\nabla g(\mathbf{z}(t))}{|\nabla g(\mathbf{z}(t))|}
\right]^2
\frac{1}{|\nabla g(\mathbf{z}(t))|}
\nabla g(\mathbf{z}(t)) \cdot \nabla r(\mathbf{z}(t))
dt
\nonumber\\
&+&
\frac{1}{\epsilon}
\int_0^{T_g}
\left|
y_g(\mathbf{z}(t))
\right|^2
\frac{1}{|\nabla g(\mathbf{z}(t))|}
\nabla g(\mathbf{z}(t)) \cdot \nabla r(\mathbf{z}(t))
dt.
\end{eqnarray*}
We have $|\mathbf{z}^\prime(t)|=|\nabla g(\mathbf{z}(t))|$
and we can write the above expression as an integral over $\partial \Omega_g$
\begin{eqnarray*}
\Gamma_r(y_g,g,r)&=&
\frac{2}{\epsilon}
\int_{\partial\Omega_g}
\left( \nabla y_g(\sigma)\cdot
\frac{\nabla g(\sigma)}{|\nabla g(\sigma)|}
\right)
\frac{\nabla r(\sigma)}{|\nabla g(\sigma)|}
\cdot\nabla y_g(\sigma)
d\sigma
\nonumber\\
&-&
\frac{2}{\epsilon}
\int_{\partial\Omega_g}
\left(
\nabla y_g(\sigma)\cdot
\frac{\nabla g(\sigma)}{|\nabla g(\sigma)|}
\right)^2
\frac{1}{|\nabla g(\sigma)|^2}
\Big(
\nabla g(\sigma) \cdot \nabla r(\sigma)
\Big)
d\sigma
\nonumber\\
&+&
\int_{\partial\Omega_g}
j\left(\sigma,\Delta y_g(\sigma)\right)
\frac{1}{|\nabla g(\sigma)|^2}
\nabla g(\sigma) \cdot \nabla r(\sigma)
d\sigma
\nonumber\\
&+&
\frac{1}{\epsilon}
\int_{\partial\Omega_g}
\left[
\nabla y_g(\sigma)\cdot
\frac{\nabla g(\sigma)}{|\nabla g(\sigma)|}
\right]^2
\frac{1}{|\nabla g(\sigma)|^2}
\nabla g(\sigma) \cdot \nabla r(\sigma)
d\sigma
\nonumber\\
&+&
\frac{1}{\epsilon}
\int_{\partial\Omega_g}
\left|
y_g(\sigma)
\right|^2
\frac{1}{|\nabla g(\sigma)|^2}
\nabla g(\sigma) \cdot \nabla r(\sigma)
d\sigma.
\end{eqnarray*}

In FreeFem++ \cite{freefem++}, for $r_h=\sum_{i=1}^{n_g} R_i \phi_i^g \in \mathbb{W}_h^g$,
with $R=(R_i)_{1\leq i \leq n_g}$,
there exist the discrete first order derivative operators
$\partial_1^h(r_h)\in \mathbb{W}_h^g$ and $\partial_2^h(r_h)\in \mathbb{W}_h^g$.
We have $\partial_1^h(r_h)=\sum_{i=1}^{n_g} R_i \partial_1^h(\phi_i^g)$, then
$\partial_1^h(r_h)$ is expressed in function of $R\in \mathbb{R}^{n_g}$.
In a similar way, we have $\partial_2^h:\mathbb{W}_h^g\rightarrow \mathbb{W}_h^g$.
We denote by $\Gamma_r^h(y_h,g_h,r_h)$ the expression obtained
from  $\Gamma_r(y_g,g,r)$
by replacing
$\partial\Omega_g,y_g,g$ by
$\partial\Omega_h, y_h,g_h$
and $\nabla r$ by $(\partial_1^h(r_h),\partial_2^h(r_h))^T$.

We set $R^r=(R_i^r)_{1\leq i\leq n_g} \in \mathbb{R}^{n_g}$ given by
\begin{equation}\label{4.7}
R^r_i=-\Gamma_r^h(y_h,g_h,\phi_i^g),\quad i=1,\dots,n_g 
\end{equation}
and we obtain
$$
\Gamma_r^h(y_h,g_h,r_h)=- \langle  R^r,R \rangle.
$$

\medskip
\textbf{Approximation of the terms containing $\mathbf{w}$}

We follow the technique from \cite{MT2022}. We can write
\begin{eqnarray}
\Gamma_w & = & \theta(g,r)\left[
j(\mathbf{x}_0,\Delta y_g(\mathbf{x}_0))
+\frac{1}{\epsilon} \left| y_g(\mathbf{x}_0) \right|^2
+\frac{1}{\epsilon} \left|
\nabla y_g(\mathbf{x}_0) \cdot
\frac{\nabla g(\mathbf{x}_0)}{|\nabla g(\mathbf{x}_0)|}
\right|^2
\right] | \nabla g(\mathbf{x}_0)|
\nonumber\\
&&+\int_0^{T_g} b_1(t) w_1(t) + b_2(t) w_2(t) dt .
\label{4.8}
\end{eqnarray}  

The ODE system (\ref{2.9})-(\ref{2.11}) is solved by the Runge-Kutta algorithm
of order 4, with initial condition $Z_0=\mathbf{x}_0$ and a step $\Delta t>0$.
We obtain $Z_k=(Z_k^1,Z_k^2)^T$ an approximation of
$\mathbf{z}_g(t_k)$, with $t_k=k\Delta t$.
The algorithm stops when $Z_m$ is ``close to'' $Z_0$.
We set the period $T_g=m\Delta t$ and we put 
$Z=(Z^1,Z^2)^T$,
where $Z^1 = (Z_k^1)_{1\leq k\leq m} \in \mathbb{R}^m$ and
$Z^2 = (Z_k^2)_{1\leq k\leq m}\in \mathbb{R}^m$.

The linear ODE system (\ref{3.10})-(\ref{3.12}) is solved
numerically by backward Euler scheme
\begin{eqnarray}
W_k &=& W_{k-1}
+\Delta t\,
A_k W_k + \Delta t\, c_k \label{4.9}\\
W_0 & = & (0,0)^T \label{4.10}
\end{eqnarray}
for $k=1, \dots, m$, $W_k=(W_k^1,W_k^2)^T$,
$c_k=\left(-\partial_2^h r_h(Z_k), \partial_1^h r_h(Z_k)\right)^T$ and
$$
A_k=\left(
\begin{array}{rr}
-\partial_{12}^h g_h(Z_k) & -\partial_{22}^h g_h(Z_k)\\
 \partial_{11}^h g_h(Z_k) &  \partial_{21}^h g_h(Z_k)
\end{array}
\right)
$$
where $\partial_i^h:\mathbb{W}_h^g\rightarrow \mathbb{W}_h^g$
and $\partial_{ij}^h:\mathbb{W}_h^g\rightarrow \mathbb{W}_h^u$, $i,j\in \{1,2\}$ are
the discrete first and second order derivative operators in FreeFem++.

From Proposition \ref{prop:2.2}, we have
$\theta(g,r)=-\frac{w_2(T_g)}{z_2^\prime (T_g)}$ if $z_2^\prime (T_g) \neq 0$.
We set 
\begin{eqnarray}
\mu_m  &=&  -\frac{\Delta t}{Z_m^2-Z_{m-1}^2}|\nabla g_h(Z_m)|
\nonumber\\
&&\times
\left[
j(Z_m,\Delta y_h(Z_m))
+\frac{1}{\epsilon} \left| y_h(Z_m) \right|^2
+\frac{1}{\epsilon} \left|
\nabla y_h(Z_m) \cdot
\frac{\nabla g_h(Z_m)}{|\nabla g_h(Z_m)|}
\right|^2  
\right]
\label{4.11}
\end{eqnarray}
and
\begin{equation}\label{4.12}
M_m =(I-\Delta t\,A_m^T)^{-1}(0,\mu_m)^T.
\end{equation}

As in \cite{MT2022}, we introduce the adjoint scheme of (\ref{4.9}), (\ref{4.10})
\begin{equation}\label{4.13}
-M_{k}=-M_{k-1} +\Delta t\, A_{k-1}^T M_{k-1} + \Delta t\,b_{k-1}
\end{equation}  
for $k=m,m-1, \dots, 1$, 
with $M_k=(M_k^1,M_k^2)^T$ and
$b_k$ is an approximation of\\
$\left(b_1(t_k),b_2(t_k)\right)^T$.
We put
$M^1 = (M_k^1)_{1\leq k\leq m}$ and $M^2 = (M_k^2)_{1\leq k\leq m}$ and it is proved in \cite{MT2022}
that $\Gamma_w $ can be approached by
$\Delta t\sum_{k=1}^{m} c_k \cdot M_k.$

For $r_h(\mathbf{x})=\sum_{i=1}^{n_g} R_i \phi_i^g(\mathbf{x})$,
with $R=(R_i)_{1\leq i \leq n_g} \in \mathbb{R}^{n_g}$, we can introduce
$\Pi_h^1=(\pi_{ij}^1)_{1\leq i,j \leq n_g}$,
$\Pi_h^2=(\pi_{ij}^2)_{1\leq i,j \leq n_g}$, the square
matrices of order $n_g$ defined by
$$
\partial_1^h(\phi_j^g)=\sum_{i=1}^{n_g} \pi_{ij}^1 \phi_i^g,\quad
\partial_2^h(\phi_j^g)=\sum_{i=1}^{n_g} \pi_{ij}^2 \phi_i^g,\quad
1\leq j \leq n_g.
$$
We obtain
$\partial_1^h r_h(\mathbf{x})
=\sum_{i=1}^{n_g} \left(\Pi_h^1 R\right)_i \phi_i^g(\mathbf{x})$ and
$\partial_2^h r_h(\mathbf{x})
=\sum_{i=1}^{n_g} \left(\Pi_h^2 R\right)_i \phi_i^g(\mathbf{x})$.
Also, we introduce
the $n_g \times m$ matrix
$$
\Phi(Z)= \left(\phi_i(Z_k)\right)_{1\leq i \leq n_g, 1\leq k\leq m}.
$$
Putting
\begin{equation}\label{4.14}
R^w=\Delta t(\Pi_h^2)^T \Phi(Z) M^1
-\Delta t(\Pi_h^1)^T\Phi(Z) M^2,\quad
R^w\in \mathbb{R}^{n_g}
\end{equation}  
we get that $\Gamma_w$ is approached by
$-\langle  R^w,R \rangle.$

\begin{remark}\label{rem:4.1}
Using (\ref{4.6}), (\ref{4.7}), (\ref{4.14}), a discrete version of (\ref{3.17}) is
\begin{equation}\label{4.15}
  d\mathcal{J}_{(G,U)}(R,V) = -\langle  R^q+R^r+R^w,R \rangle
  -\langle V^q,V\rangle
\end{equation}
where $R^q,R^r,R^w,R \in \mathbb{R}^{n_g}$ and $V^q,V \in \mathbb{R}^{n_u}$.
\end{remark}

\begin{proposition}\label{prop:4.2}
We assume that $R^q+R^r+R^w\neq 0$ or $V^q\neq 0$. We have
\begin{itemize}  
\item[i)] $R^*=R^q+R^r+R^w$ and $V^*= V^q$ is a descent direction for (\ref{4.15}).
\item[ii)] $R^*=K^{-1}(R^q+R^r+R^w)$ and $V^*= V^q$
  is a descent direction for (\ref{4.15}), where $K$ is a symmetric
  positive definite matrix.
\end{itemize}  
\end{proposition}

\noindent
\textbf{Proof.} i) We obtain
$d\mathcal{J}_{(G,U)}(R^*,V^*)=-\|R^*\|^2 -\|V^*\|^2 <0$.

ii) In this case $KR^*=R^q+R^r+R^w$, then
$$
d\mathcal{J}_{(G,U)}(R^*,V^*)=-\langle  KR^*,R^* \rangle -\|V^*\|^2 <0
$$
since $K$ is positive definite.
\quad$\Box$

\section{Numerical results}
\setcounter{equation}{0}
Our approach allows both shape and topology optimization. The examples involve both opening and closing of holes.
We use the software FreeFem++, \cite{freefem++}.
We set $D=]-3,3[\times ]-3,3[$,
$j\left(\sigma,\Delta y_g(\sigma)\right)
=\frac{1}{2}\left(\Delta y_g(\sigma)\right)^2$.
The finite elements for $g_h$ is 
$\mathbb{P}_3$, for $u_h$ is $\mathbb{P}_1$ and
for $y_h$ we use the Hsieh, Clough and Tocher (HCT) finite element, \cite{Dautray1990}.

At each iteration $k$, for the line-search, we evaluate the cost function 30 times or more,
for $\lambda=\rho^i\lambda_0$, $i=0,\dots 29$ with $\rho=0.8$ and $\lambda_0>0$ and
we choose the global minimum obtained for these values of $\lambda$. This type of backtracking line-search does not use the Armijo rule, \cite{Bertsekas1999}, and allows several local minimum points in the computed values.
The tolerance parameter for the stopping test is $tol$.

\medskip

\textbf{Test 1.}

\begin{figure}[ht]
\begin{center}
  \includegraphics[width=7cm]{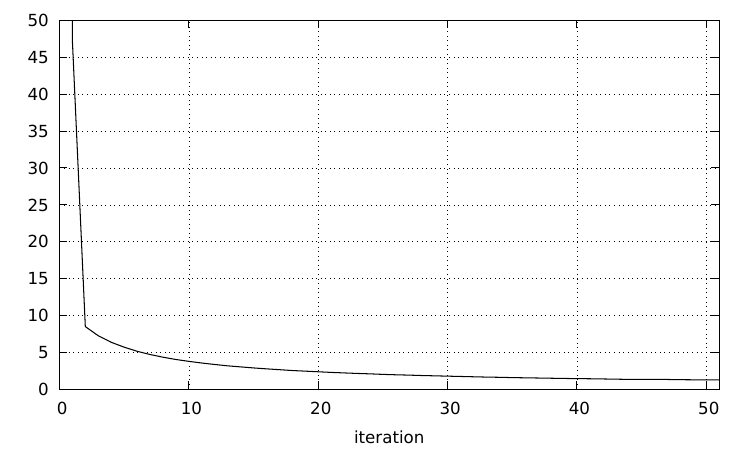}
\end{center}
\caption{Test 1. The penalized cost function for iterations $k\geq 1$.
  At $k=0$ the cost function is $1310.71$ and the last value is $1.25341$.}
\label{fig:test2_J}
\end{figure}

\begin{figure}[ht]
\begin{center}
  \includegraphics[width=6cm]{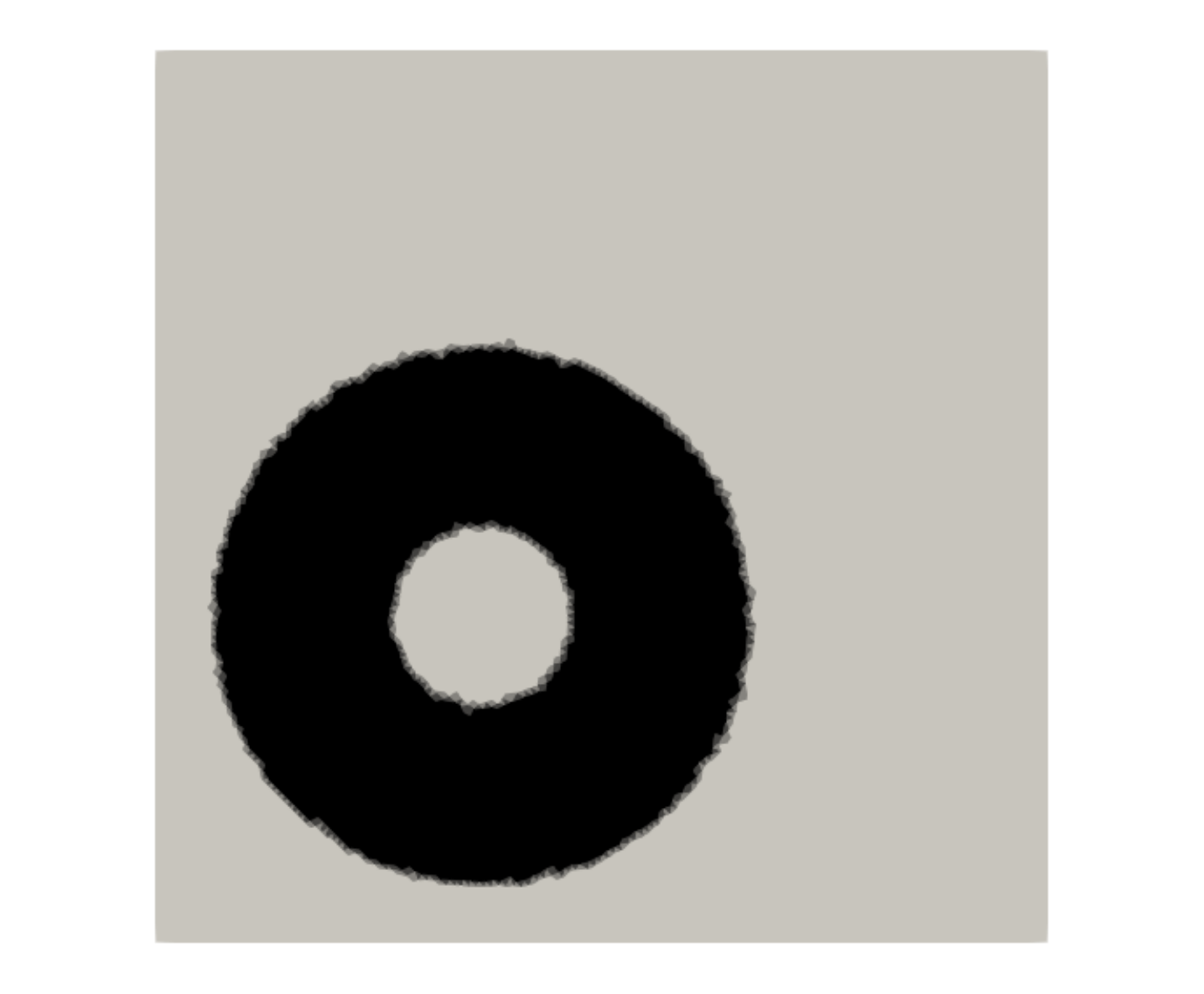}
  \includegraphics[width=6cm]{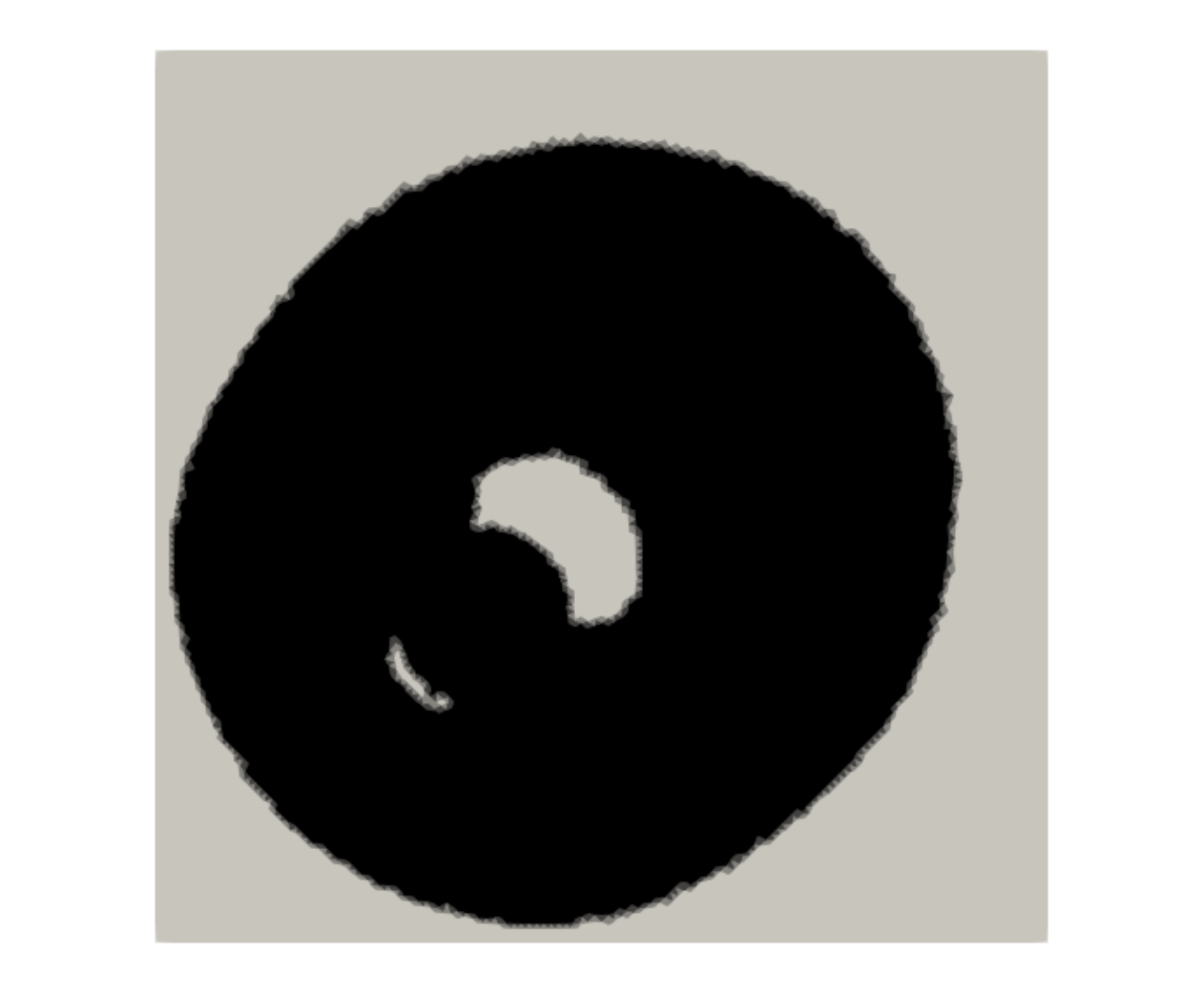}\\
  \includegraphics[width=6cm]{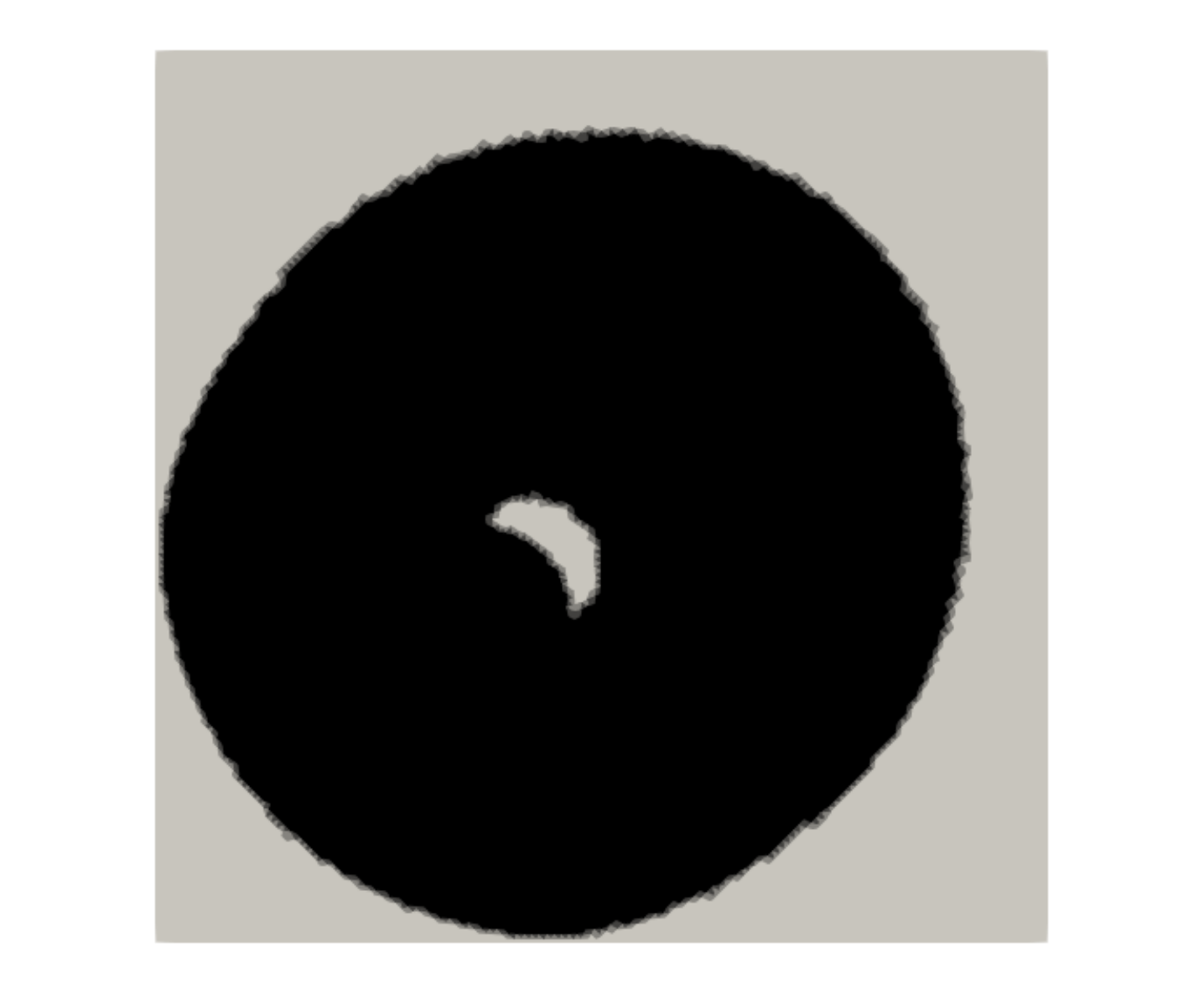}
  \includegraphics[width=6cm]{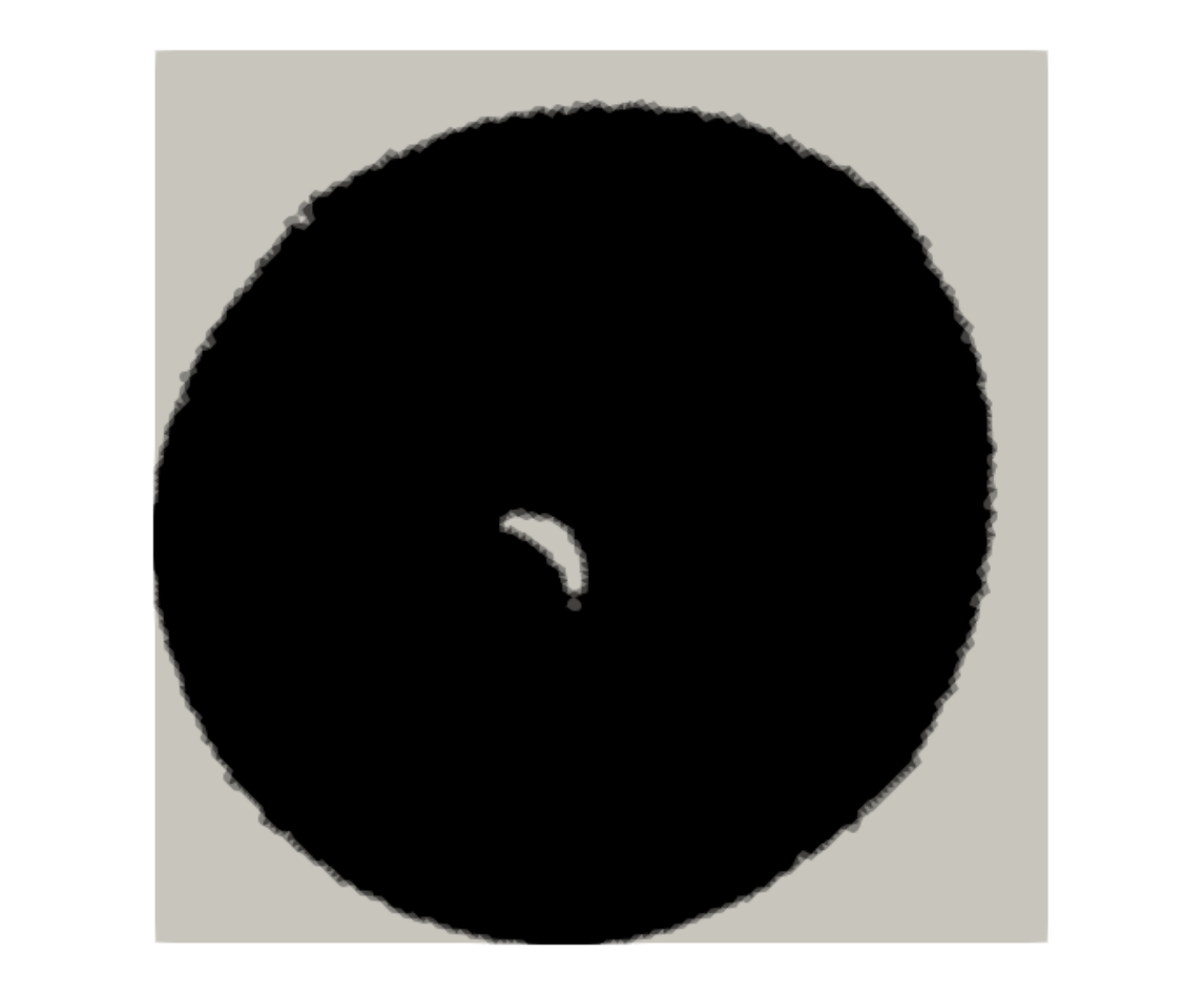}
\end{center}
\caption{Test 1. The domains $\Omega_k$ for k=0, 1, 2, final.}
\label{fig:test2_Omega}
\end{figure}

We set $f=0$. Then, by (\ref{2.1}), (\ref{2.2}) we see that for any admissible choice of the domain $\Omega$
 we have $y_{\Omega} = 0$ everywhere, that is the cost (\ref{2.15}) is constant, i.e. any admissible $\Omega$
is optimal. This example is specially conceived to test the validity of our finite element routine.
The penalization parameter is $\epsilon=0.8$ (see (\ref{3.7})) and
the initial domain $\Omega_0$ is given by
$$
g_0(x_1,x_2)=
\max\left(
-(x_1+0.8)^2-(x_2+0.8)^2+0.6^2,
(x_1+0.8)^2 +(x_2+0.8)^2-1.8^2
\right), 
$$
a disk of center $(-0.8,-0.8)$ and radius $1.8$
with a circular hole of center $(-0.8,-0.8)$ and radius $0.6$.
The initial guess for the control is $u_0=1$ and notice that for this choice of
$g_0, u_0$ the discrete penalized cost (\ref{4.1})  has the value
1310.71. The algorithm performs a very strong decrease, see Fig. \ref{fig:test2_J},
and in the first iteration it creates a supplementary hole
(that is closed afterwards), see Fig. \ref{fig:test2_Omega}.
Consistent shape modifications are also obtained.
The triangulation has 17175 vertices and 33868 triangles.
For the line-search we use $\lambda_0=1$ and $tol=10^{-6}$.

We use a simplified direction $(R^k,V^k)$ given by
$$
R^k=R^{q,k}+R^{r,k},\quad V^k=V^{q,k}
$$
where $R^{q,k}$, $V^{q,k}$, $R^{r,k}$ are obtained using the formulas
(\ref{4.6}), (\ref{4.7}), with $(G^k,U^k)$ in place of $(G,U)$.
This simplified direction is obtained from to the case i) of Proposition
\ref{prop:4.2}, by neglecting the term $R^w$ and
it can be implemented easier in FreeFem++.
There is neither regularization nor normalization of the simplified descent direction.
The integrals on the boundary of
$\Omega_g$ can be computed with the command 
\texttt{int1d(Th,levelset=gh)(\dots)}.
The algorithm stops after 51 iterations.

%\clearpage

In Figure \ref{fig:test2_J} we show the history of the penalized cost function who
has three terms $\mathcal{J}=t_1 +\frac{1}{\epsilon}t_2 +\frac{1}{\epsilon}t_3$.
At the final iteration, we have:
$t_1=0.825763$, $t_2=0.226998$, $t_3=0.115119$, $\mathcal{J}_{51}=1.25341$.
We recall that $t_1$ is the original cost function and $t_2$, $t_3$ ``small'' means that
the boundary condition (\ref{2.2}) holds.

\begin{figure}[ht]
\begin{center}
  \includegraphics[width=9cm]{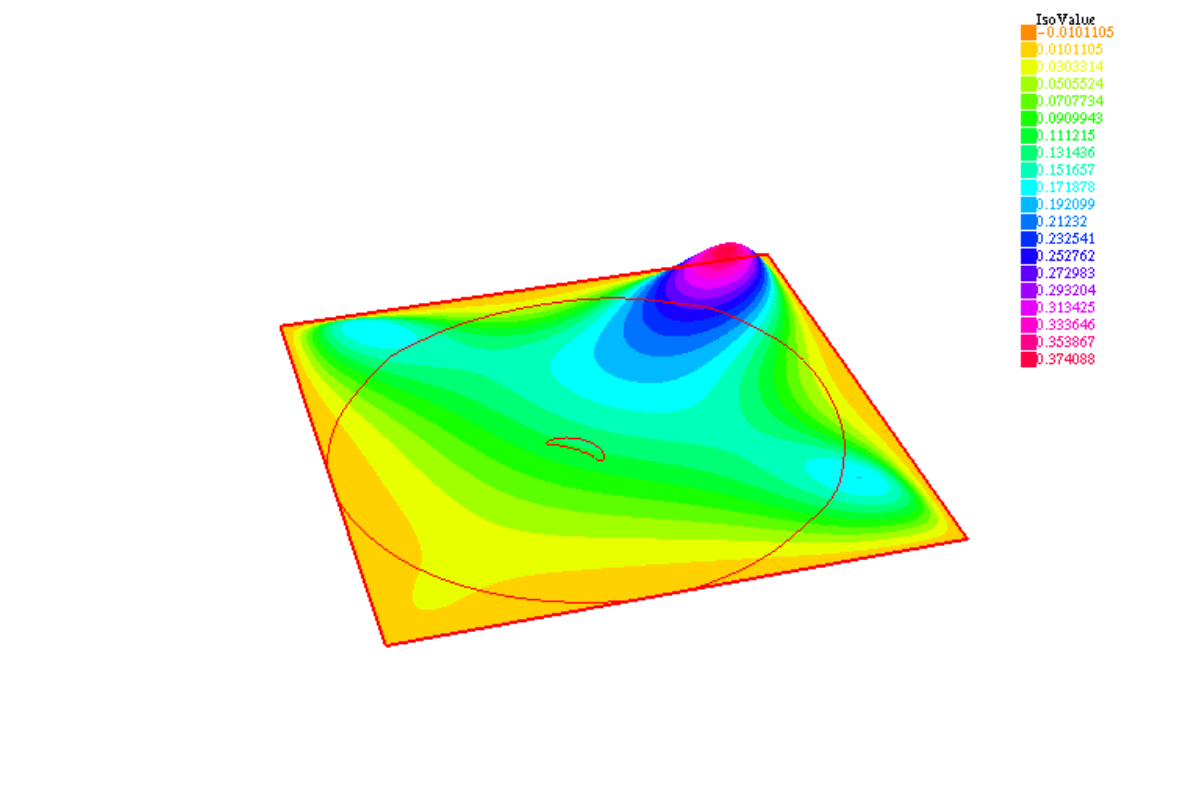}\\
  \includegraphics[width=6.5cm]{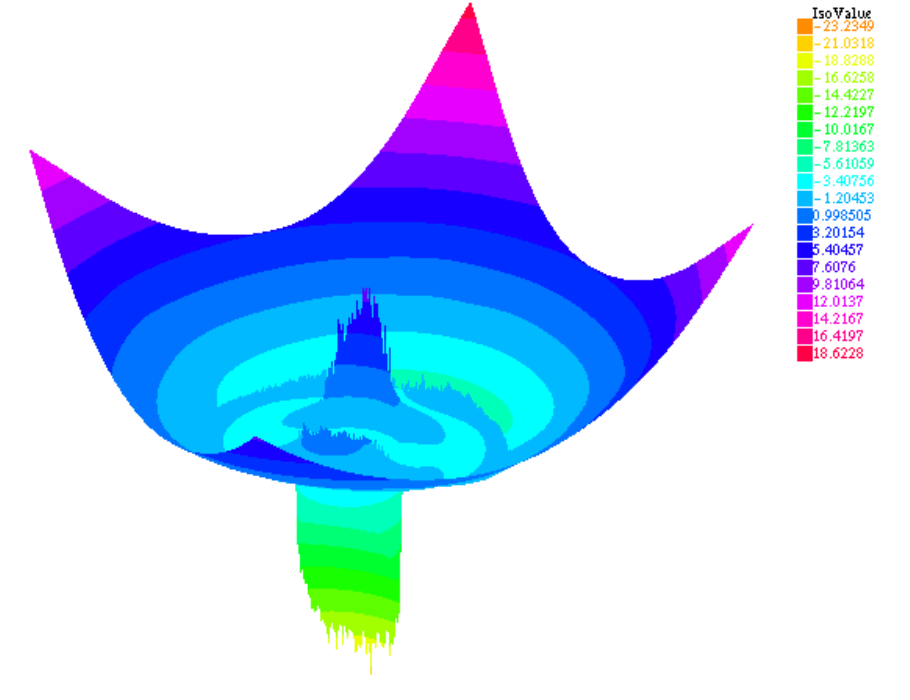}
  \includegraphics[width=6.5cm]{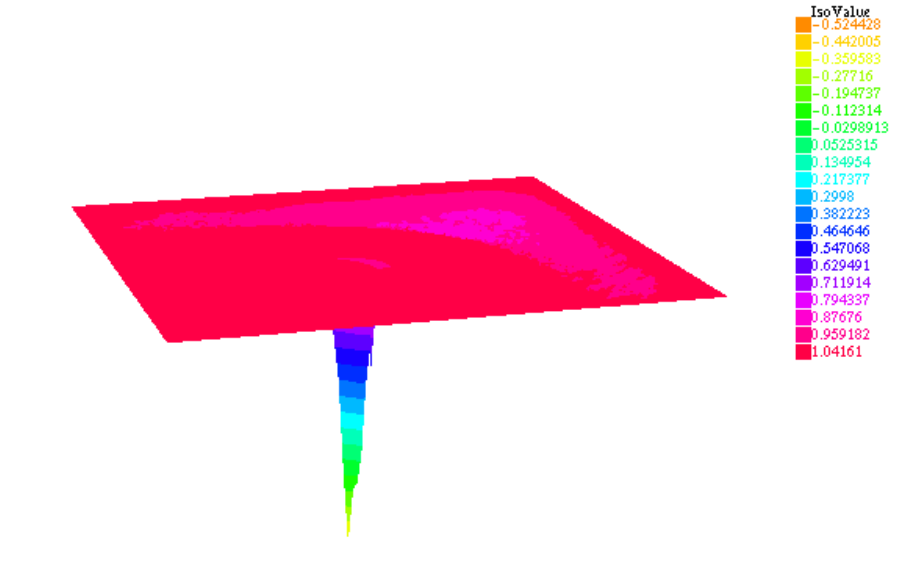}
\end{center}
\caption{Test 1. Final $y_h$ and $\partial\Omega_h$ (top).
Final $g_h$ (left, bottom) and $u_h$ (right, bottom).}
\label{fig:test2_y}
\end{figure}

The final $y_h$, $g_h$, $u_h$ are presented in Figure \ref{fig:test2_y}.
The parametrization function $g_h:D\rightarrow\mathbb{R}$ is globally continuous,
but we observe that it is not smooth enough.
In the next test, we also employ a smoothness procedure.

\clearpage

\textbf{Test 2.}

\textbf{Case a)}

This case is as before, in particular $\epsilon=0.8$, but we change other parameters:  the direction $R^k$,
the triangulation has 9698 vertices and 19034 triangles, $f=0.1$,
$\lambda_0=0.5$ and $tol=10^{-2}$.

Let $\overline{r}_h^k$ be the $\mathbb{P}_3$ finite elements function
associated to the vector $\overline{R}^k=R^{q,k}+R^{r,k}+R^{w,k}$,
where $R^{q,k}$, $R^{r,k}$ $R^{w,k}$ are
obtained using the formulas (\ref{4.6}), (\ref{4.7}), (\ref{4.14})
with $(G^k,U^k)$ in place of $(G,U)$.
We point out that we use the whole gradient for the descent direction.
We have observed that $\overline{r}_h^k$ is not smooth enough.
Then, inspired by \cite{Burger2003}, we solve the problem:
find $\tilde{r}_h^k\in \mathbb{V}_h$,
such that
$$
\int_D \Delta \tilde{r}_h^k  \Delta \varphi_h \, d\mathbf{x}
= \int_D  \overline{r}_h^k \varphi_h \, d\mathbf{x},
\qquad \forall \varphi_h \in \mathbb{V}_h
$$
and we set $R^k$ as the vector associated to $\Pi_{\mathbb{P}_3} (\tilde{r}_h^k)$,
the $\mathbb{P}_3$ projection of $\tilde{r}_h^k$.

This direction is inspired by the case ii) of Proposition
\ref{prop:4.2}. If we neglect the effect of interpolation between $\mathbb{P}_3$ and
the HCT finite element, we have
$$
R^k=K^{-1}M\overline{R}^k=K^{-1}M(R^{q,k}+R^{r,k}+R^{w,k})
$$
where $K$, $M$ are two symmetric, positive definite matrices.
Using HCT in place of $\mathbb{P}_3$ finite element at the regularization stage
improves the smoothness of the domain parametrization function $g_h$.

\begin{figure}[ht]
\begin{center}
\includegraphics[width=7cm]{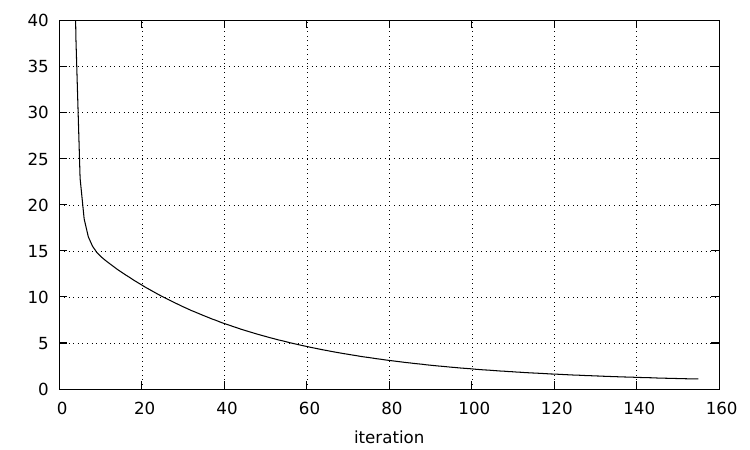}
\end{center}
\caption{Test 2, case a). The penalized cost function for iterations $k\geq 4$.
  The first values are: $\mathcal{J}_{0}=1344.08$, $\mathcal{J}_{1}=724.423$,
  $\mathcal{J}_{2}=184.156$, $\mathcal{J}_{3}=74.8713$ and the last value is $1.09214$.}
\label{fig:test4_J}
\end{figure}

\begin{figure}[ht]
\begin{center}
  \includegraphics[width=6cm]{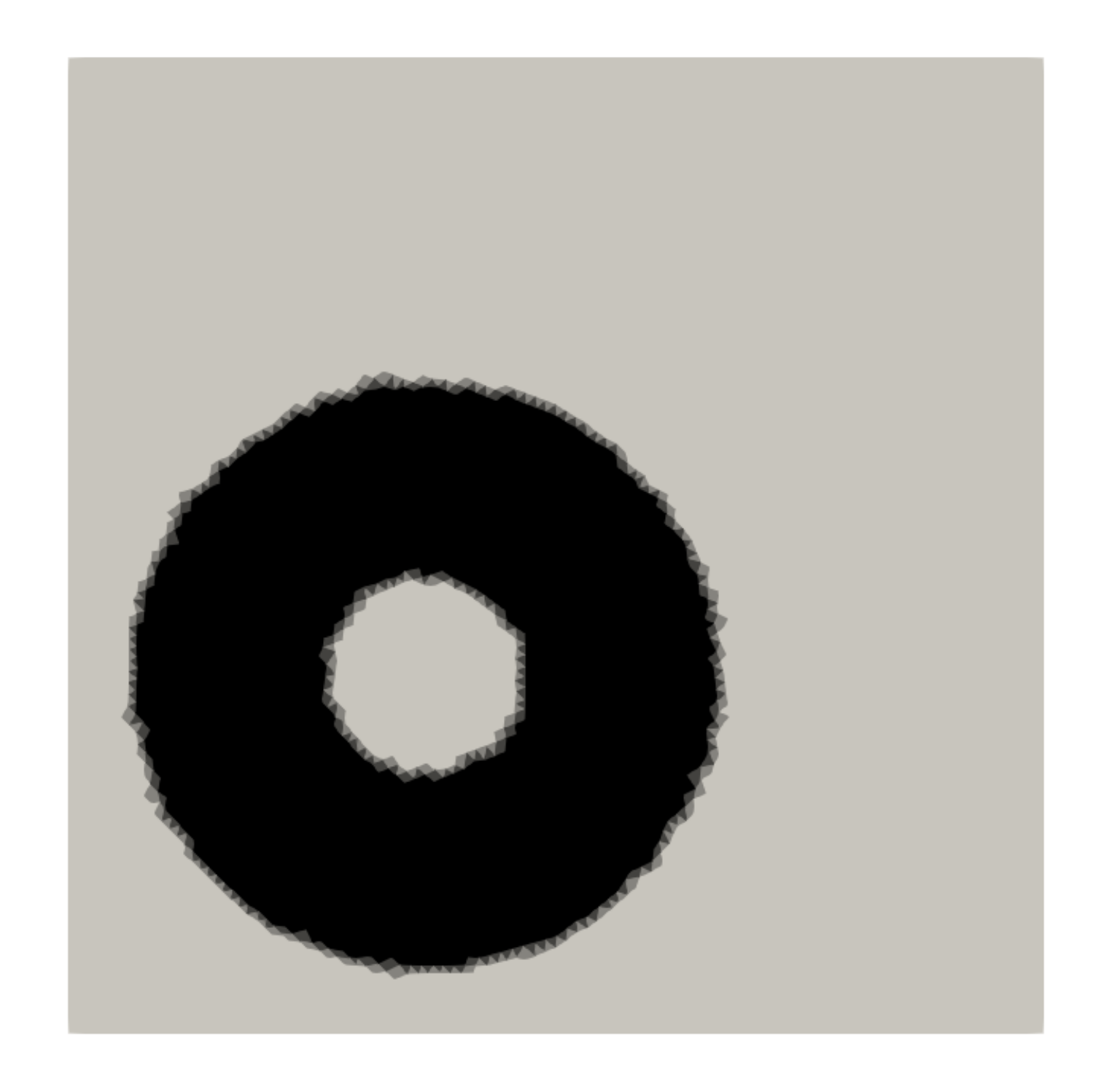}
  \includegraphics[width=6cm]{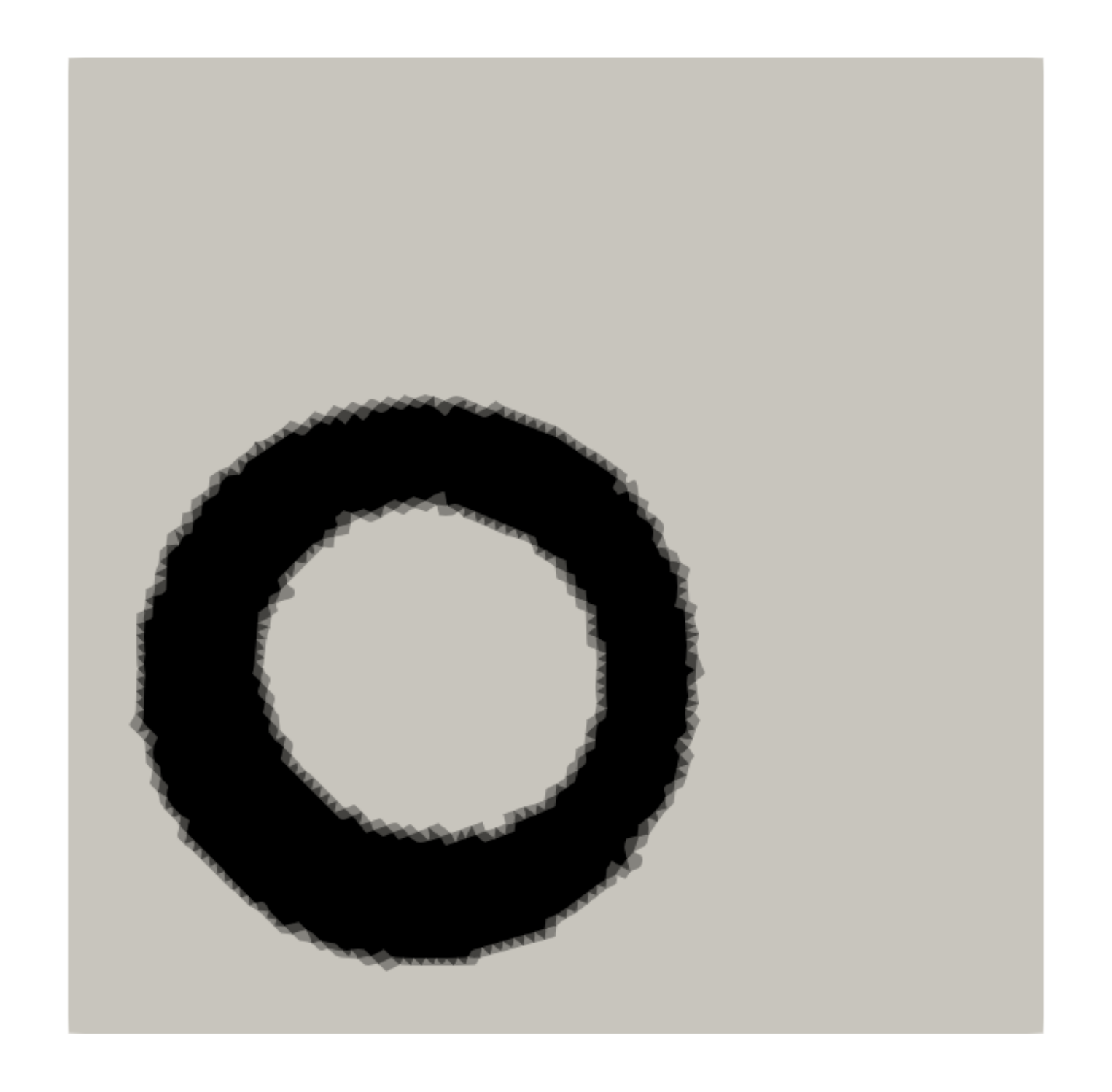}\\
  \includegraphics[width=6cm]{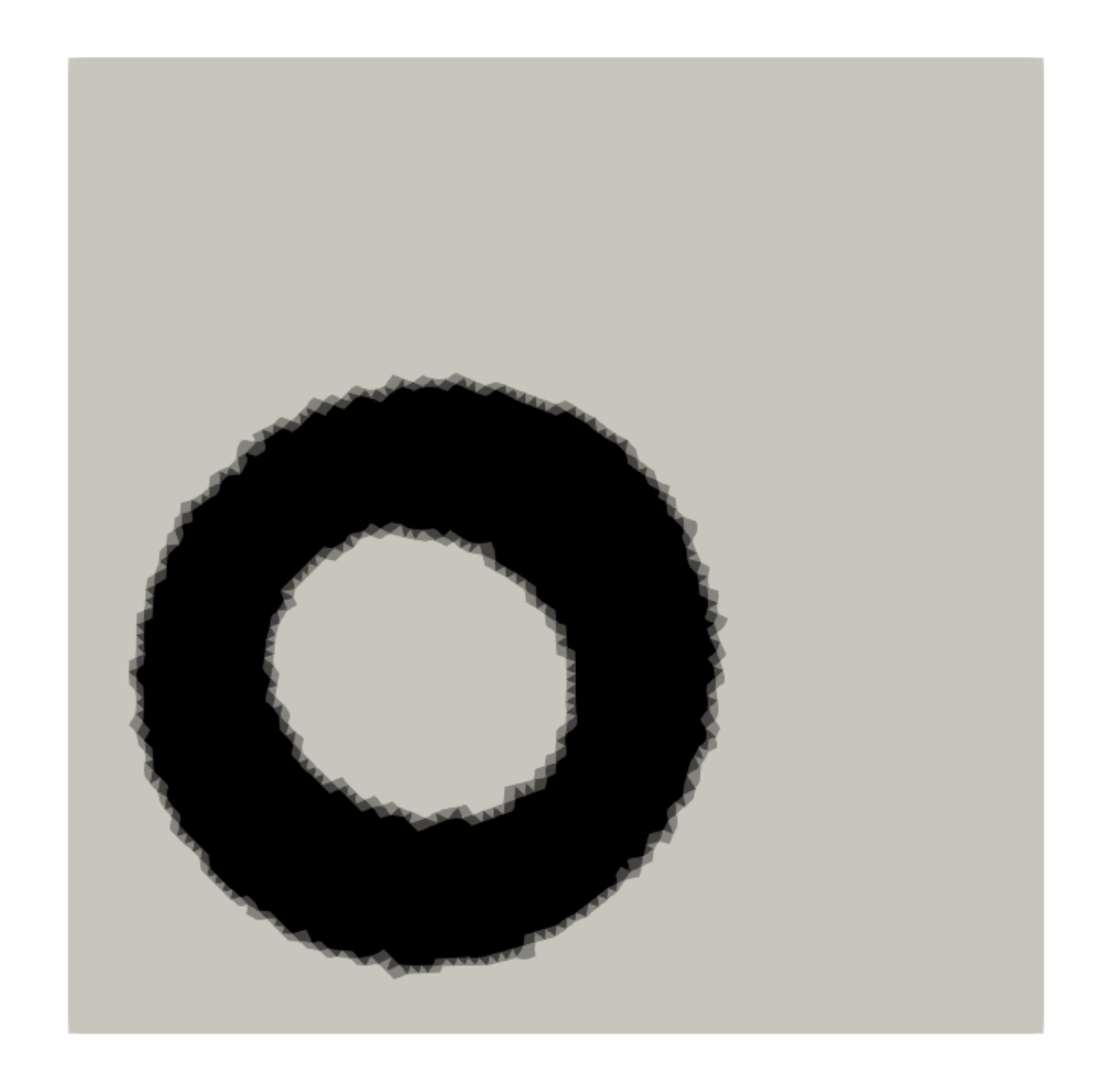}
  \includegraphics[width=6cm]{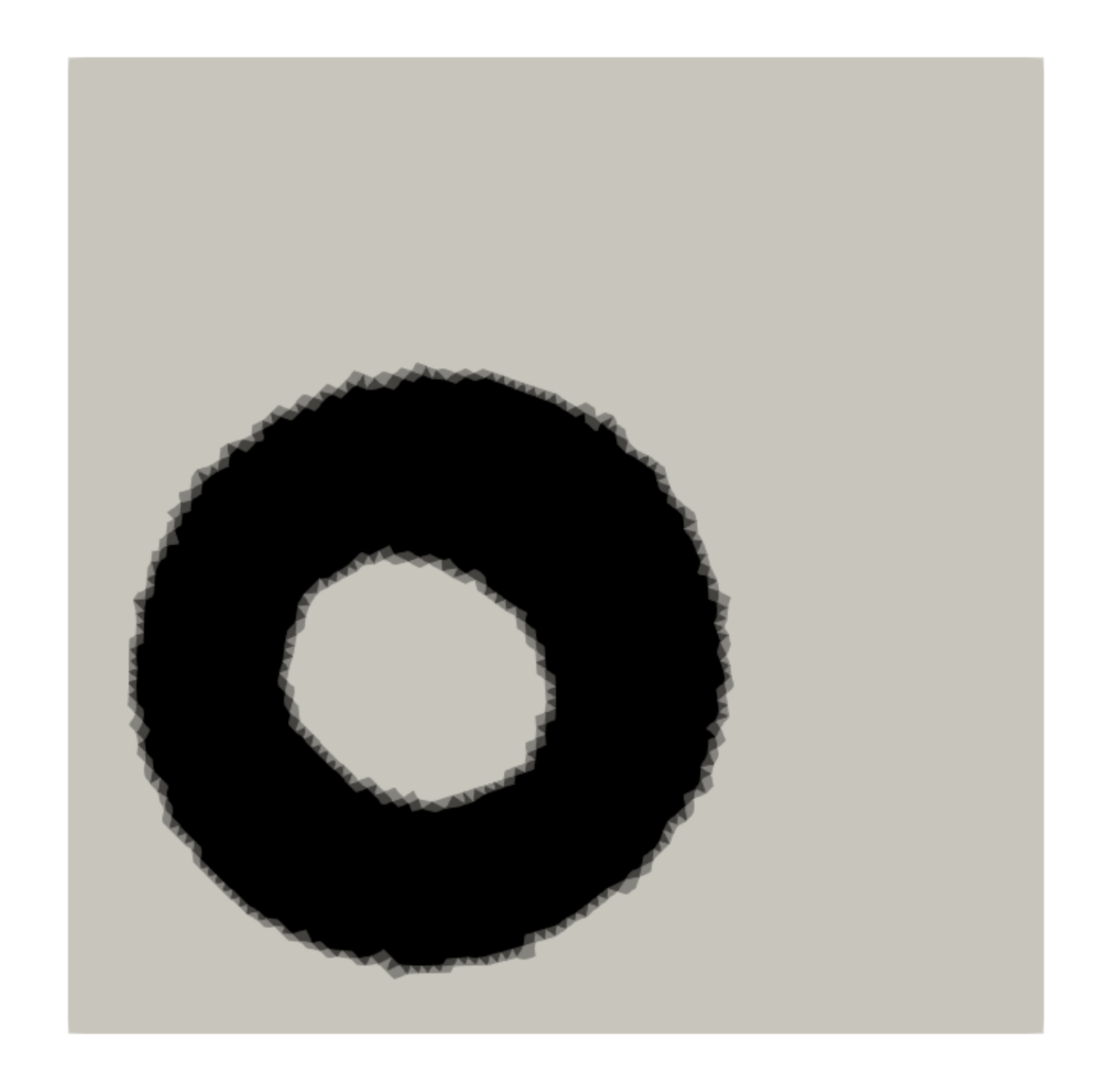}
\end{center}
\caption{Test 2, case a). The domains $\Omega_k$ for k=0, 1, 2, 155.}
\label{fig:test4_Omega}
\end{figure}

The algorithm stops after 155 iterations.
We observe in Figure \ref{fig:test4_J} the history of the penalized cost function
and in Figure \ref{fig:test4_Omega} some $\Omega_k$.
After 155 iterations we have: $t_1=0.496559$, $t_2=0.289739$, $t_3=0.18673$ and
$\mathcal{J}_{155}=t_1 +\frac{1}{\epsilon}t_2 +\frac{1}{\epsilon}t_3=1.09214$.

\begin{figure}[ht]
\begin{center}
  \includegraphics[width=9cm]{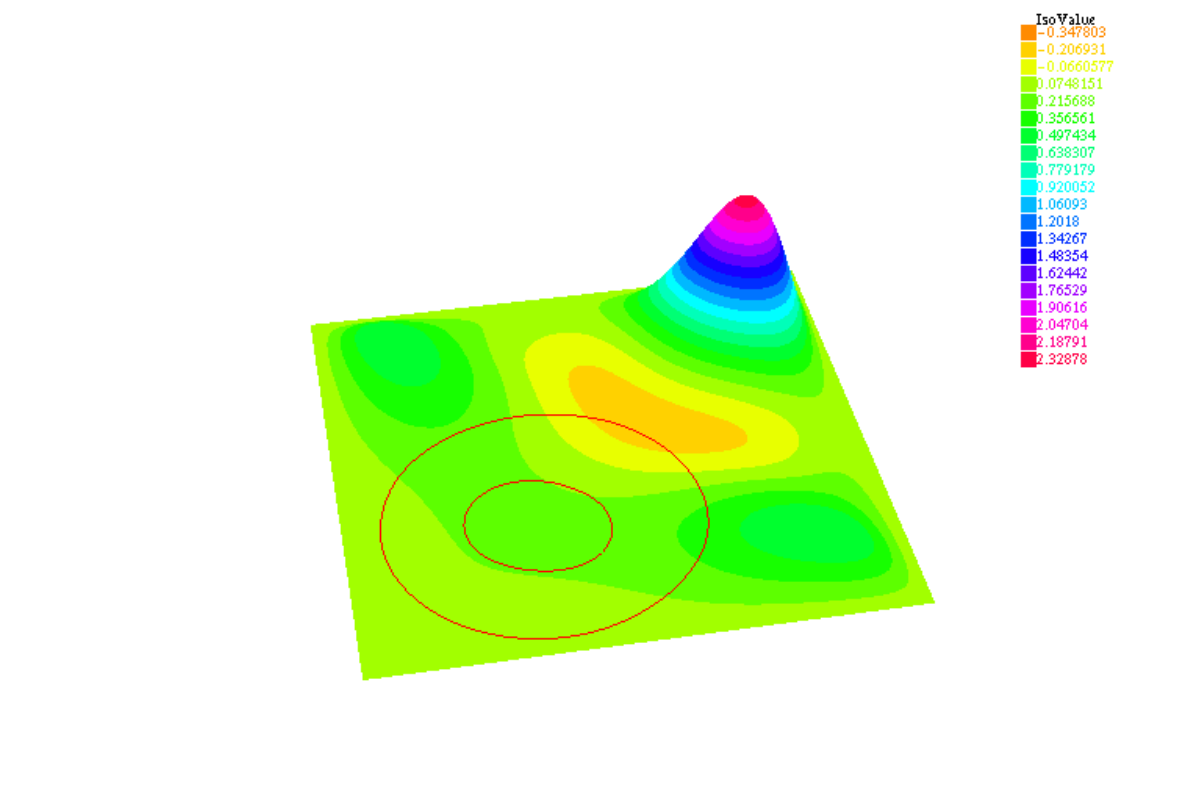}\\
  \includegraphics[width=6.5cm]{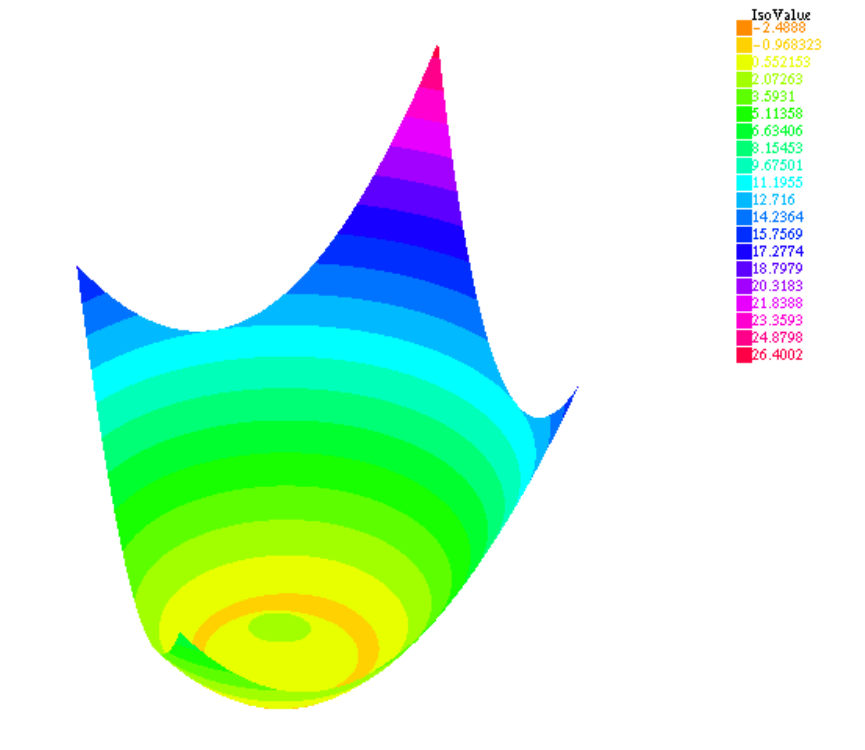}
  \includegraphics[width=6.5cm]{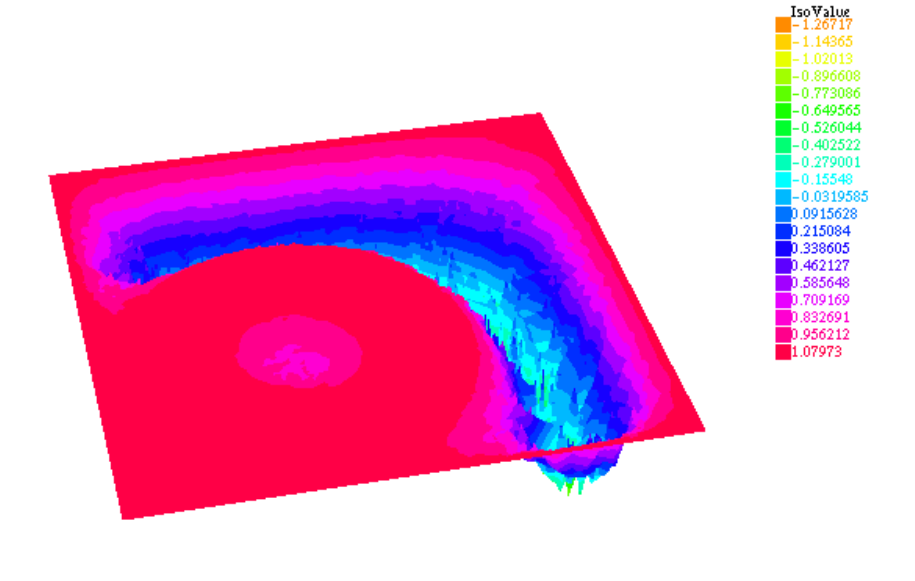}
\end{center}
\caption{Test 2, case a). Final $y_h$ and $\partial\Omega_h$ (top).
  Final $g_h$ (left, bottom) and $u_h$ (right, bottom).}
\label{fig:test4_y}
\end{figure}

In Figure \ref{fig:test4_y} we plotted the final $y_h$, $g_h$, $u_h$.
This time, $g_h$ is smoother than before. The matrices $K^{-1}$ and $M$ operate
as filter to improve the smoothness of $\overline{r}_h^k$.
As before, $u_h$ is non constant in the exterior of $\Omega_h$.

\begin{figure}[ht]
\begin{center}
  \includegraphics[width=9cm]{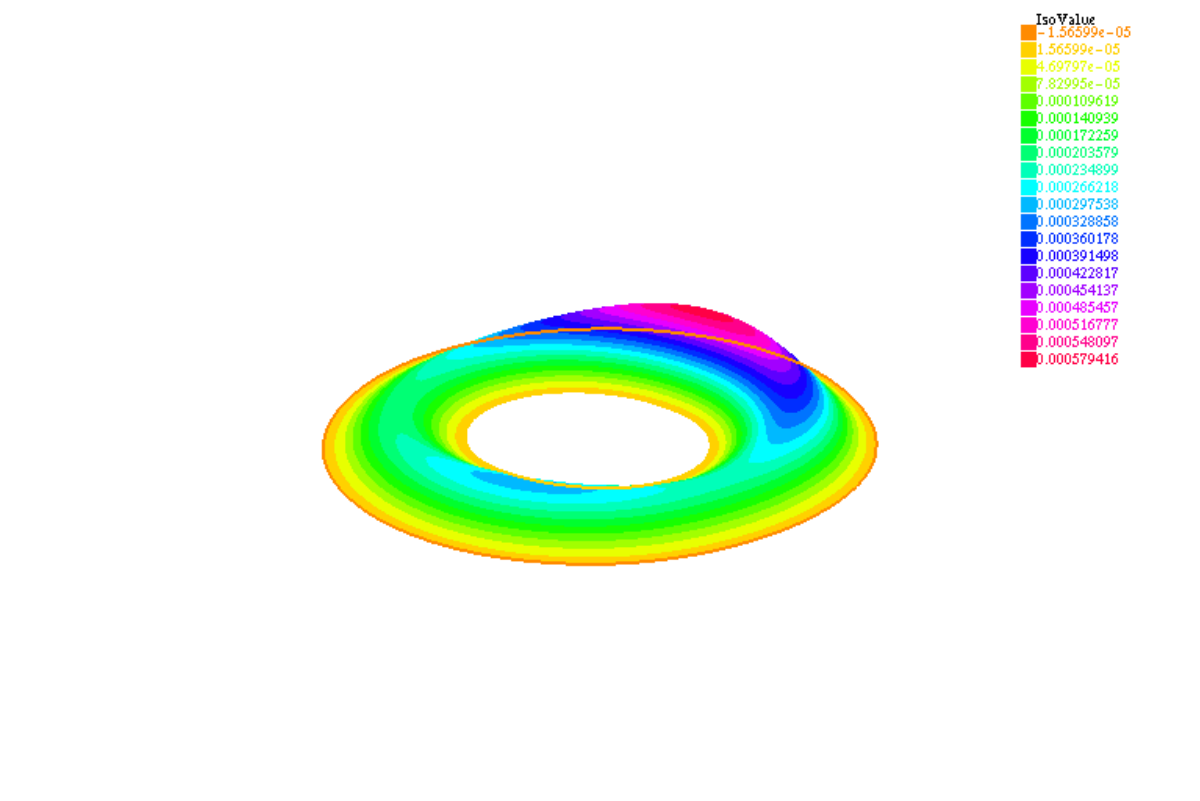}
\end{center}
\caption{Test 2, case a). Solution of the clamped plate model in the final domain.}
\label{fig:test4_y1}
\end{figure}

We have also computed the solution of (\ref{2.1})-(\ref{2.2}) for the final domain,
see Figure \ref{fig:test4_y1}.

\clearpage

\textbf{Case b)}

We use $\epsilon=0.01$, $\lambda_0=1$ and
the initial domain $\Omega_0$, corresponding to $g_0(x_1,x_2)$ which is
maximum of the following
functions: $(x_1)^2 + (x_2)^2 -2.5^2$, $-(x_1+1)^2 - (x_2+1)^2 +0.6^2$,
$-(x_1-1)^2 - (x_2+1)^2 +0.6^2$, $-(x_1+1)^2 - (x_2-1)^2 +0.6^2$
(the disk of center $(0,0)$ and radius $2.5$
with three circular holes of radius $0.6$).
The other values are as in the Test2, case a), in particular
the descent direction is given by the whole gradient and we use the same regularization of this
direction. Moreover, we normalize the descent direction
as follows
$$
\left(
\frac{\overline{R}^k}{\|\overline{R}^k\|_\infty},
\frac{V^k}{\|V^k\|_\infty}
\right).
$$

\begin{figure}[ht]
\begin{center}
\includegraphics[width=7cm]{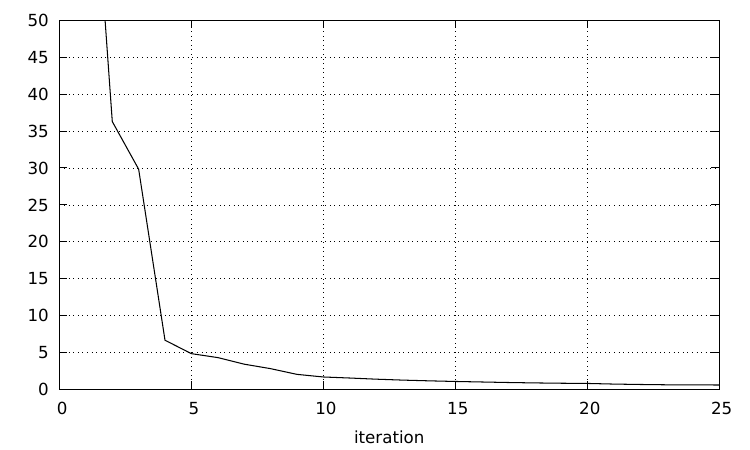}
\end{center}
\caption{Test 2, case b). The penalized cost function for iterations $k\geq 2$.
  The first values are: $\mathcal{J}_{0}=257.585$, $\mathcal{J}_{1}=86.8112$
and the last value is $0.544788$.}
\label{fig:test2_b_J}
\end{figure}

The algorithm stops after 25 iterations. The history of the penalized cost function
is presented in Figure \ref{fig:test2_b_J}.

\begin{figure}[ht]
\begin{center}
  \includegraphics[width=4cm]{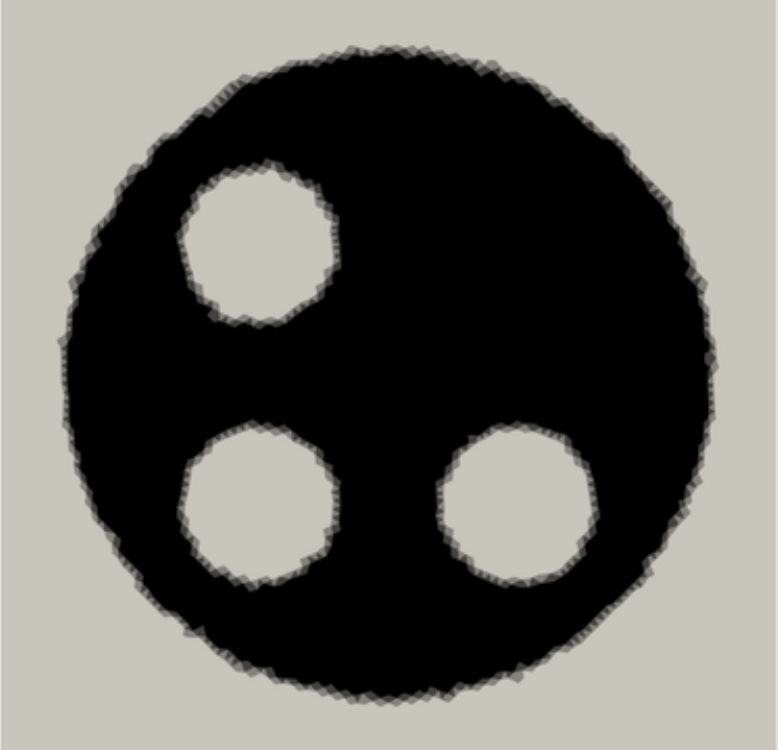}
  \includegraphics[width=4cm]{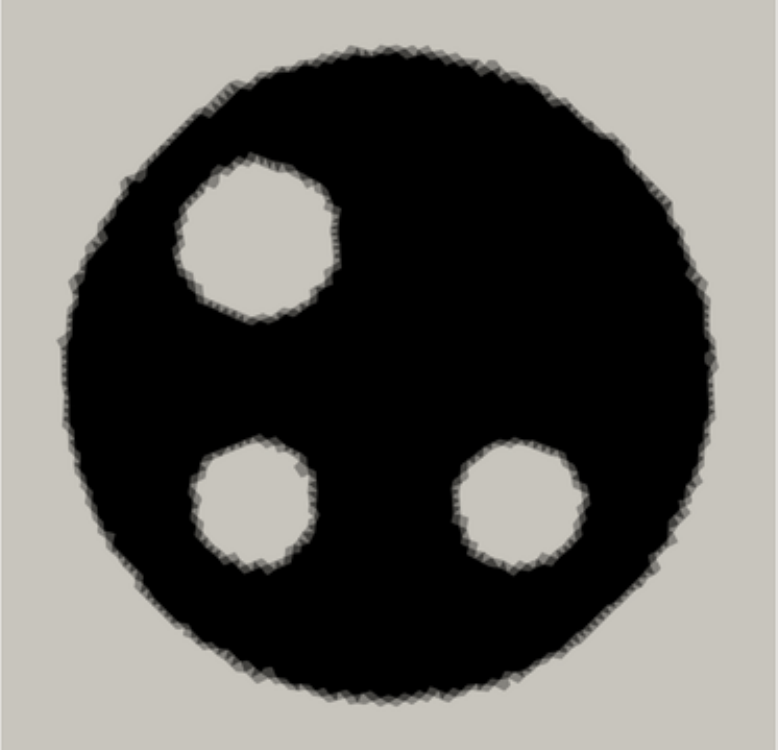}
  \includegraphics[width=4cm]{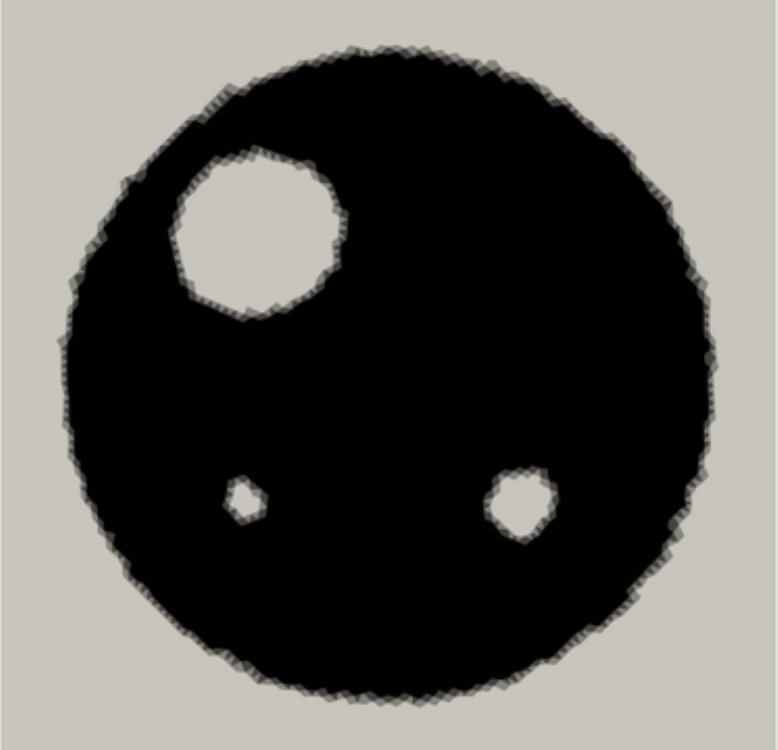}\\
  \includegraphics[width=4cm]{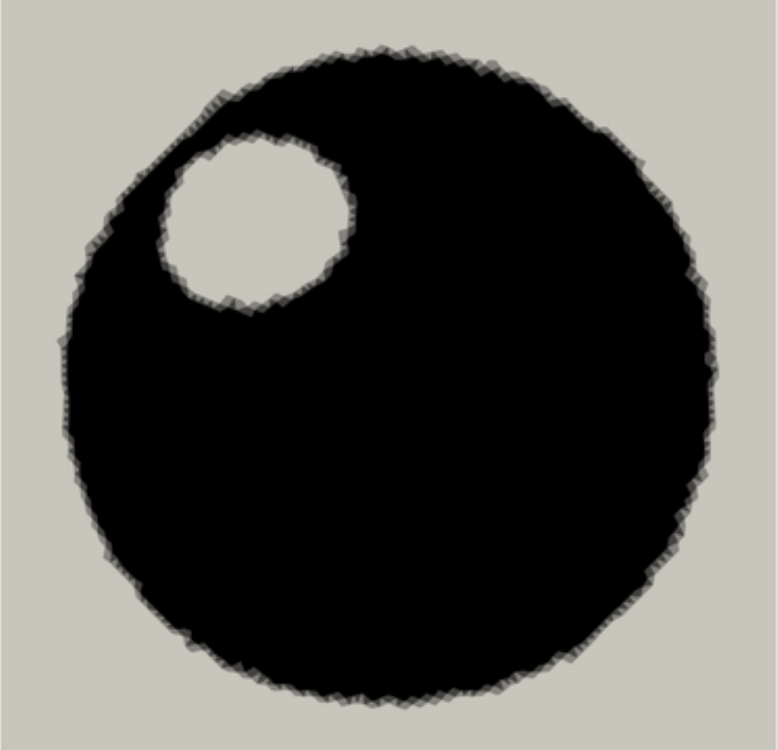}
  \includegraphics[width=4cm]{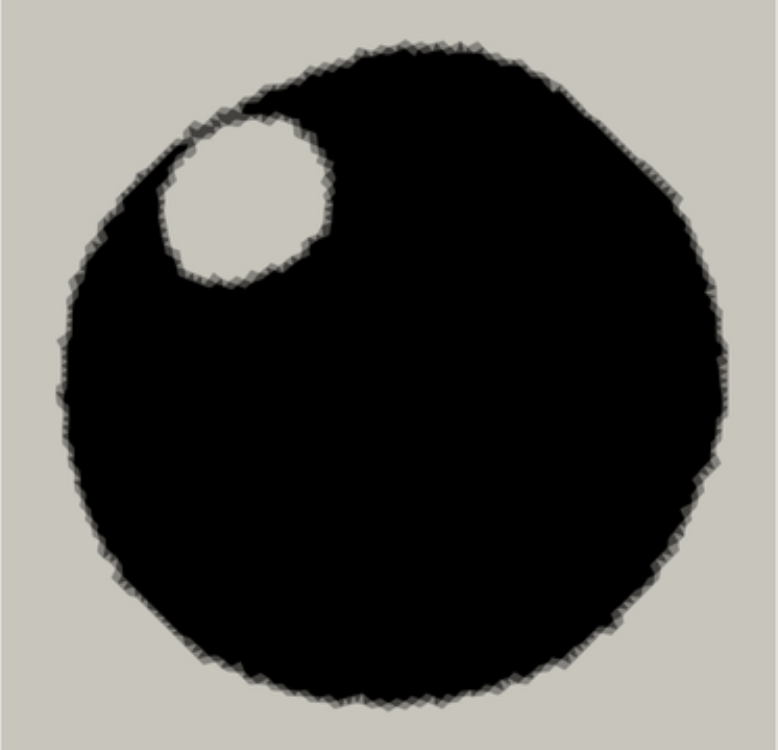}
  \includegraphics[width=4cm]{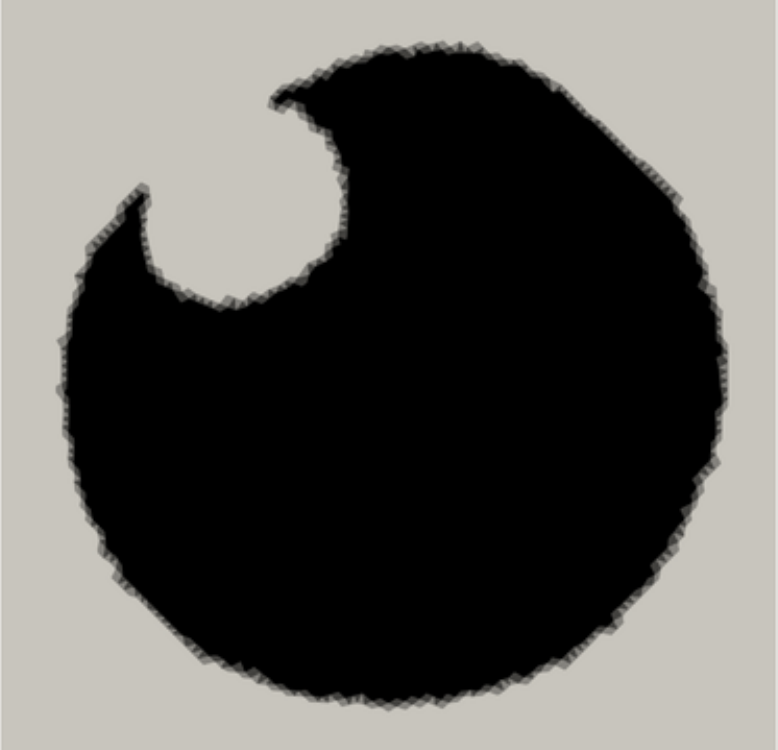}\\
  \includegraphics[width=4cm]{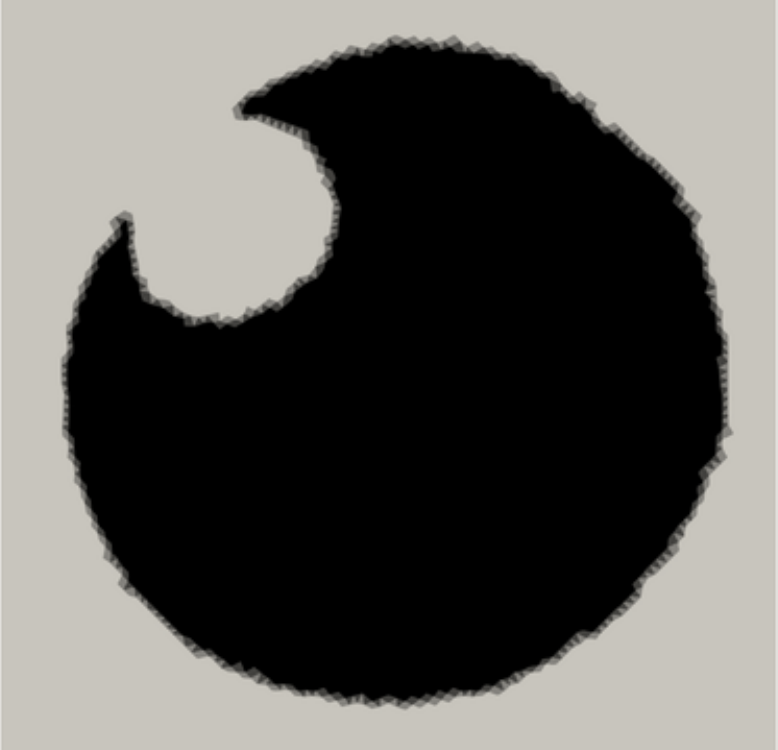}
  \includegraphics[width=4cm]{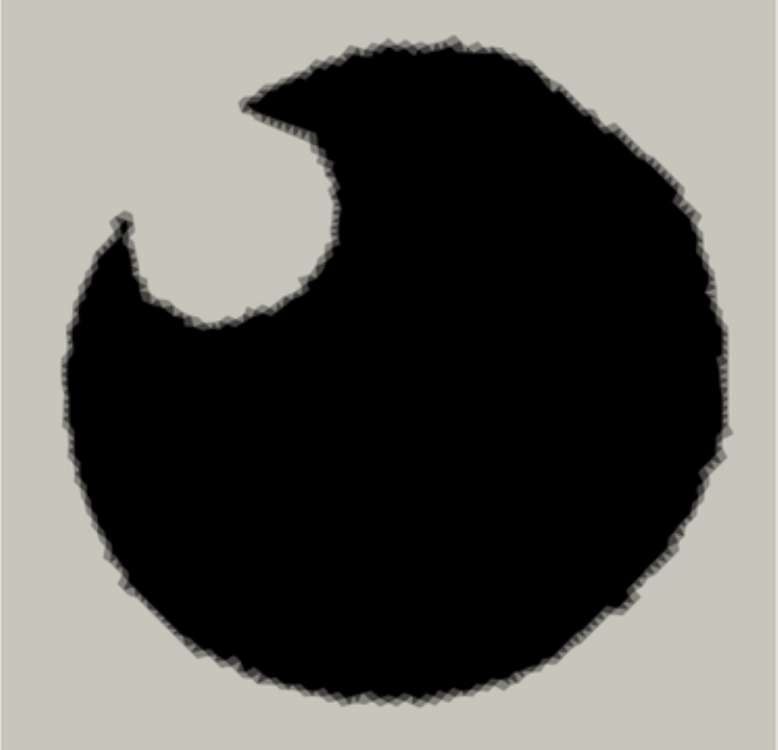}
  \includegraphics[width=4cm]{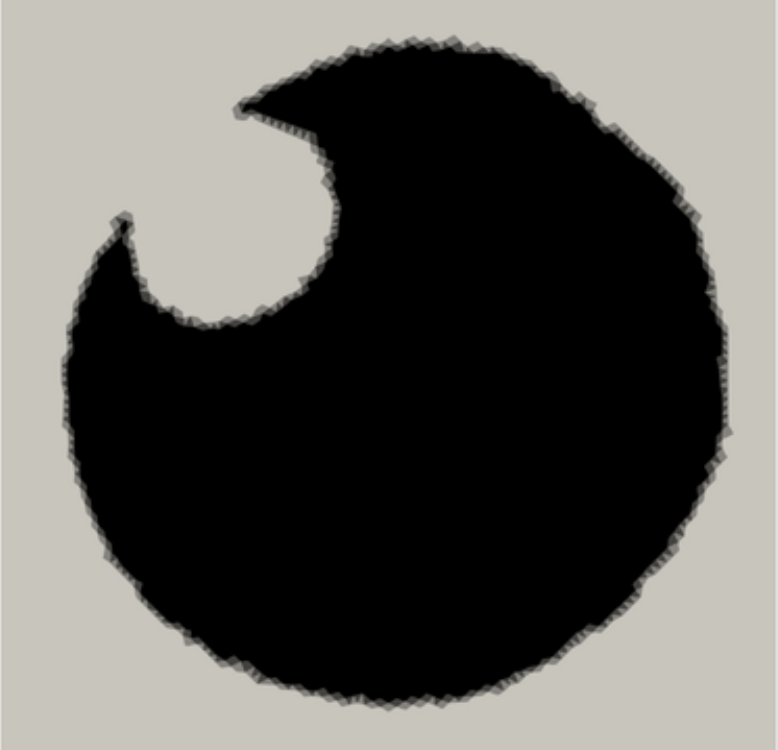}
\end{center}
\caption{Test 2, case b). The domains $\Omega_k$ for $k=0$ (top, left),
  some intermediate domains during the line-search after $k=0$ ($\lambda=0.8^8$, top, middle
  $\lambda=0.8^4$, top, right),
  $k=1,2,3$ (middle line),
  $k=10,17,25$ (bottom line).}
\label{fig:test2_b_Omega}
\end{figure}

We plotted some $\Omega_k$ in Figure \ref{fig:test2_b_Omega}, initially there are 3 holes,
at $k=1$ only one hole and finally the domain is simply connected.

In Table \ref{tab:test2_b}, we show the terms $t_1$, $t_2$, $t_3$ composing the
penalized cost function. We observe that the values on the line $t_1$ are decreasing.
The terms $t_2$, $t_3$ concerning the boundary conditions (\ref{2.2}) are small at
the final iteration.

\begin{figure}[ht]
  \begin{center}
  \includegraphics[width=7cm]{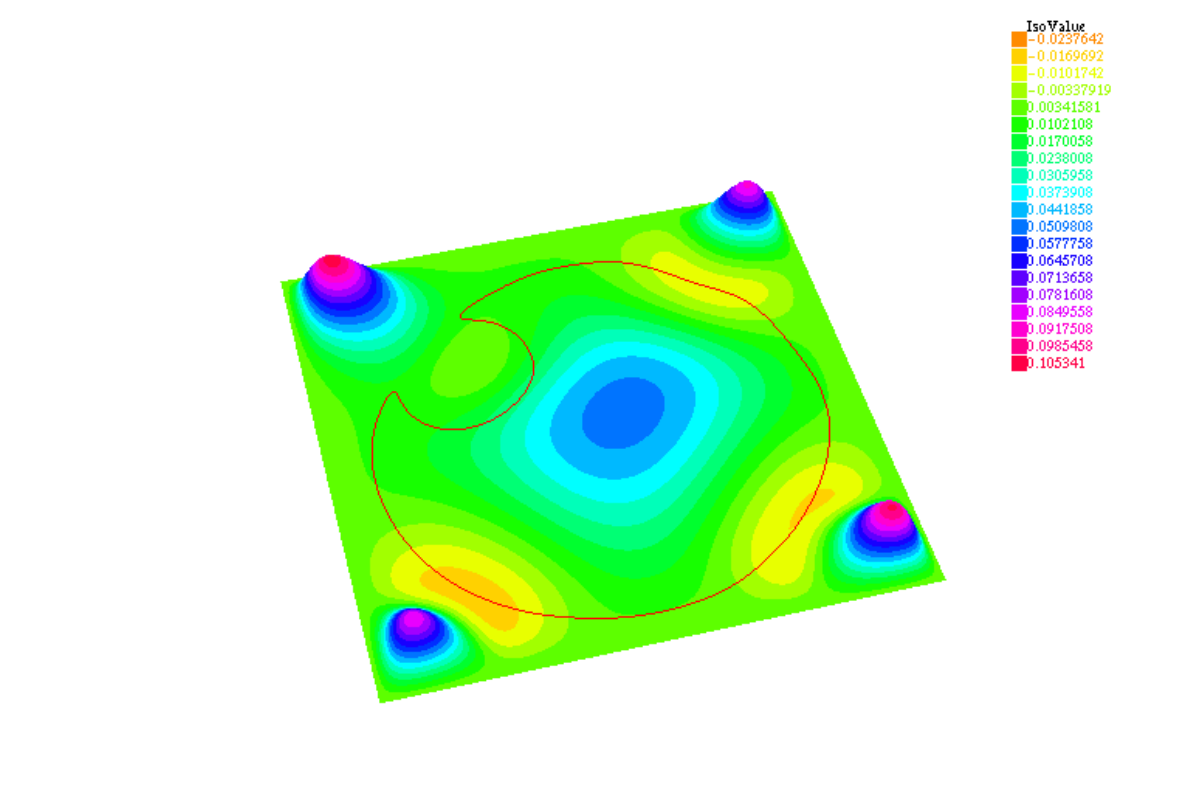}  
  \includegraphics[width=7cm]{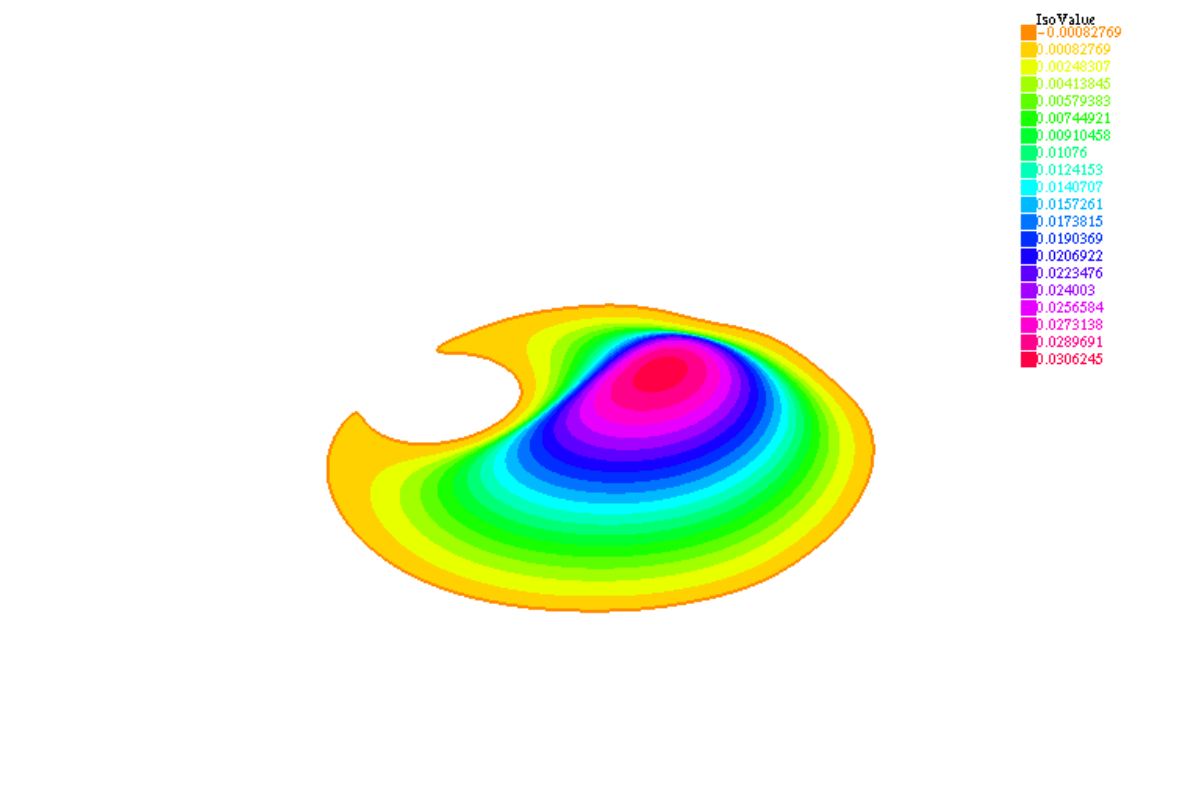}
\end{center}
  \caption{Test 2, case b). Final $y_h$ and $\partial \Omega_h$ (left) and
    the solution of the clamped plate model in the final domain (right).}
\label{fig:test2_b_yh_y1}
\end{figure}

We have computed the solution of (\ref{2.1})-(\ref{2.2}) for the final domain,
see Figure \ref{fig:test2_b_yh_y1}.

\begin{table}[ht]
\begin{center}
\begin{tabular}{|c|r|r|r|r|r|r|}
\hline
it.            & k=0     &         &         & k=1     & k=2     &   k=3   \\  \hline
$t_1$          & 1.06091 & 0.90691 & 0.69644 & 0.38435 & 0.12072 & 0.07963 \\  \hline
$t_2$          & 1.89564 & 1.539   & 0.99339 & 0.50898 & 0.17096 & 0.14157 \\  \hline
$t_3$          & 0.66960 & 0.57684 & 0.48017 & 0.35528 & 0.19079 & 0.15602 \\  \hline
$\mathcal{J}$  & 257.585 & 212.491 & 148.054 & 86.8112 & 36.2963 & 29.8393 \\  \hline
\end{tabular}

\begin{tabular}{|c|r|r|r|}
\hline
it.            & k=10    &   k=17  & k=25    \\  \hline
$t_1$          & 0.05954 & 0.05354 & 0.05215  \\  \hline
$t_2$          & 0.00478 & 0.00315 & 0.00182  \\  \hline
$t_3$          & 0.01103 & 0.00514 & 0.00309  \\  \hline
$\mathcal{J}$  & 1.64094 & 0.88375 & 0.54478  \\  \hline
\end{tabular}
\end{center}
\caption{Test 2, case b). The computed objective function 
  $\mathcal{J}=t_1 +\frac{1}{\epsilon}t_2 +\frac{1}{\epsilon}t_3$ for the domains
  as in Figure \ref{fig:test2_b_Omega}. 
}
\label{tab:test2_b}
\end{table}

\clearpage

\medskip

\textbf{Test 3.}

Now, the penalization parameter is $\epsilon=0.1$, $f=0.1$.
The initial domain is as in Test 2, case b).
The initial guess for the control is $u_0=1$ and $tol=10^{-6}$. 
We use a triangulation of 17175 vertices and 33868 triangles.

\begin{figure}[ht]
\begin{center}
\includegraphics[width=7cm]{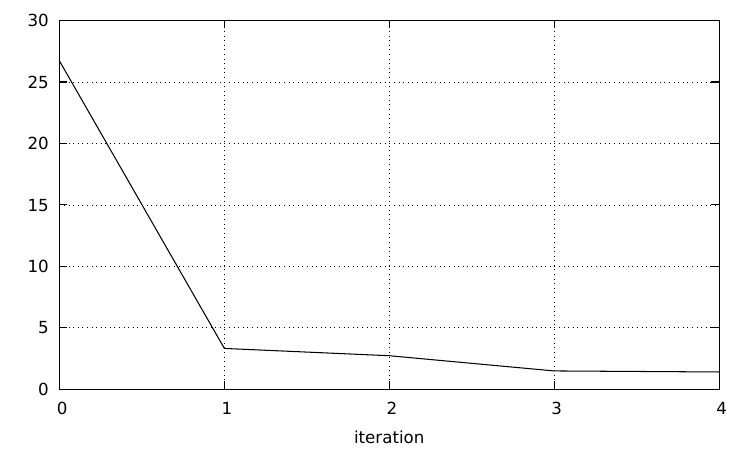}
\end{center}
\caption{Test 3. The penalized cost function decreases from $26.7086$ to $1.39281$.}
\label{fig:test7_J}
\end{figure}

For this test, we use a simplified direction $(R^k,V^k)$ corresponding
to
\begin{equation}\label{5.1}
r_h^k=\sum_{i=1}^{n_g} R_i^k \phi_i^g = \Pi_{\mathbb{P}_3} (-p_h u_h),\quad
v_h^k=\sum_{i=1}^{n_u} V_i^k \phi_i^u= \Pi_{\mathbb{P}_1} (-p_h)
\end{equation}
the $\mathbb{P}_3$ interpolation of $-p_h u_h$ and 
the $\mathbb{P}_1$ interpolation of $-p_h$, respectively.
It is different from (\ref{4.6}), but when the quadrature formula employs
only the values of functions in the vertices, we have
\begin{eqnarray*}
&&\int_D \left( 2(g_h)_+ u_h p_h r_h^k + (g_h)_+^2 p_h v_h^k \right)
d\mathbf{x}
=-\int_D  2(g_h)_+ u_h p_h \Pi_{\mathbb{P}_3} (p_h u_h)d\mathbf{x}\\
&&-\int_D (g_h)_+^2 p_h \Pi_{\mathbb{P}_1} (p_h) d\mathbf{x}
=-\int_D  2(g_h)_+ \left( \Pi_{\mathbb{P}_3} (p_h u_h) \right)^2 d\mathbf{x}
-\int_D (g_h)_+^2 \left( \Pi_{\mathbb{P}_1} (p_h) \right)^2 d\mathbf{x} \leq 0.
\end{eqnarray*}
From (\ref{4.5}) and the before inequality, we obtain that $\Gamma_q$, the
terms containing $q$ in (\ref{3.17}), is negative.
In this numerical test, we have also noted that this simplified choice is
even a descent direction.
There is neither regularization nor normalization of the simplified descent direction.
The algorithm stops after 4 iterations.
See Figure \ref{fig:test7_J}, with the history of the penalized cost function.
At the final iteration, we have: $t_1=0.287516$, $t_2=0.0631412$, $t_3=0.047388$ and
$\mathcal{J}_4=t_1 +\frac{1}{\epsilon}t_2 +\frac{1}{\epsilon}t_3=1.39281$.

\begin{figure}[ht]
\begin{center}
  \includegraphics[width=4cm]{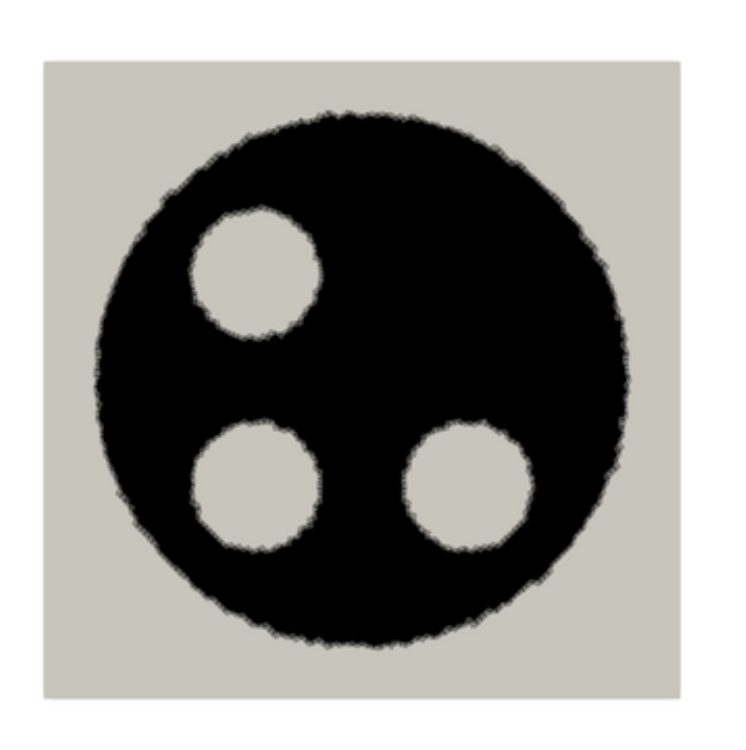}
  \includegraphics[width=4cm]{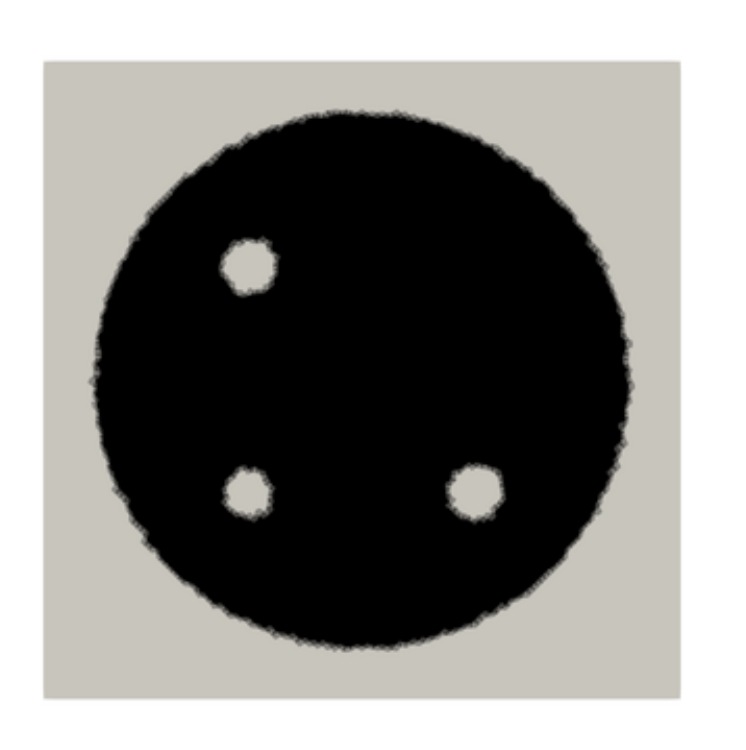}
  \includegraphics[width=4cm]{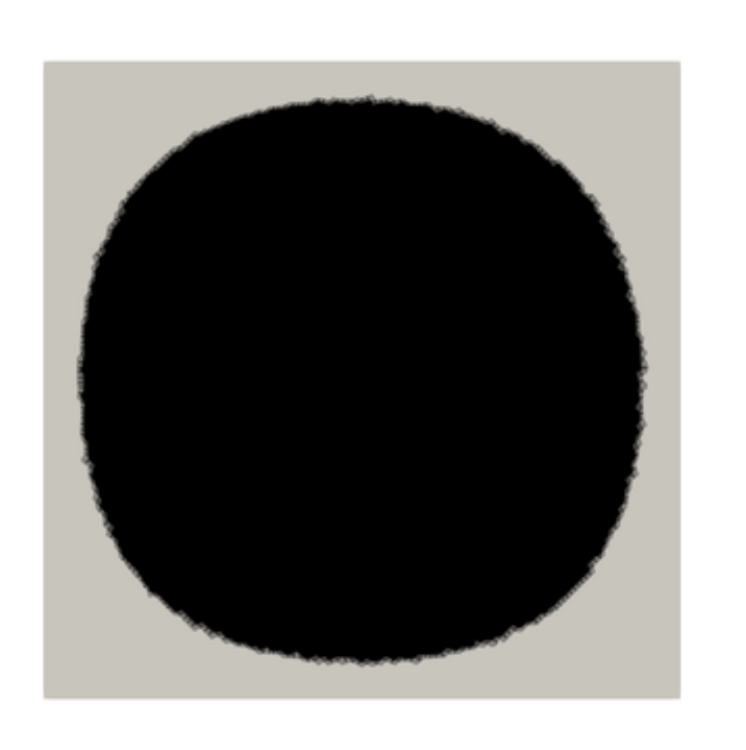}\\
  \includegraphics[width=4cm]{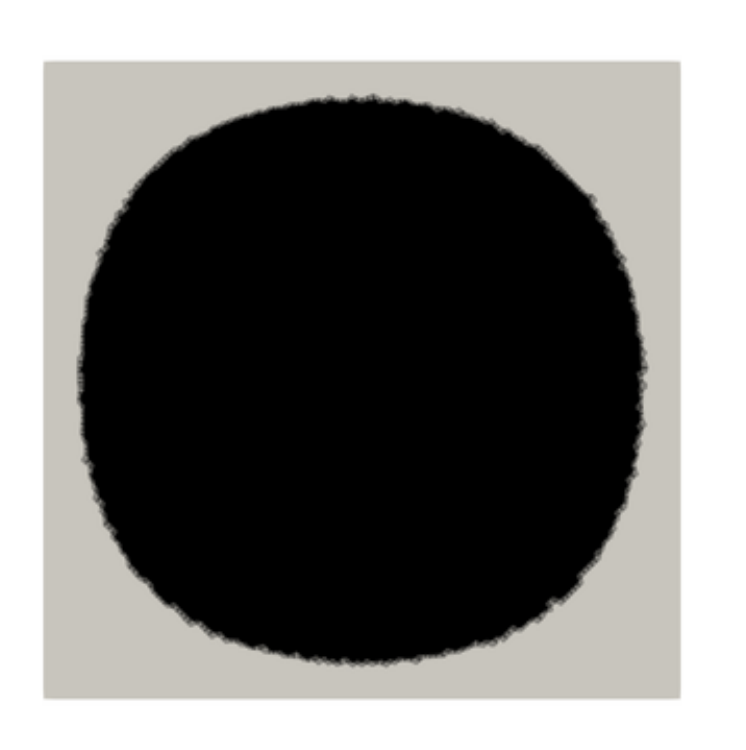}
  \includegraphics[width=4cm]{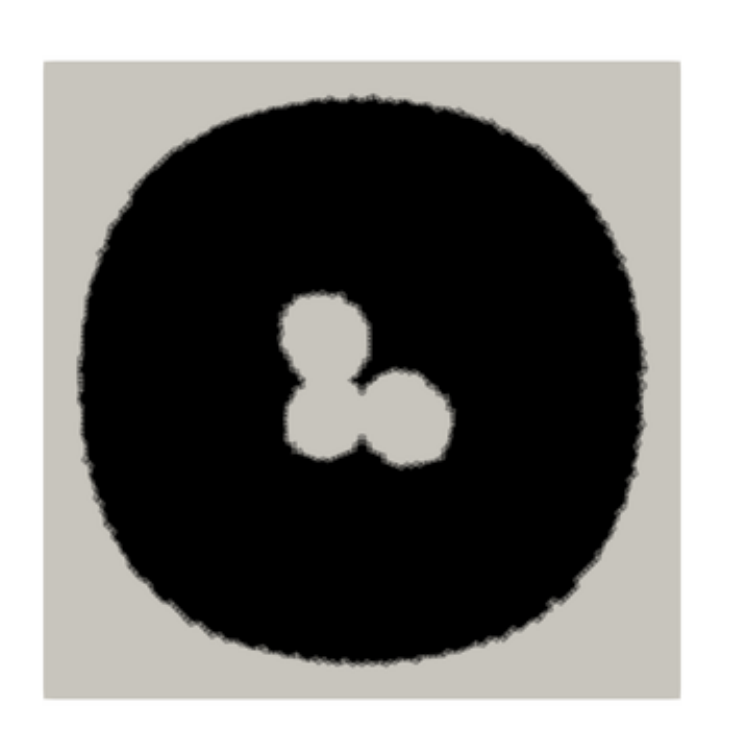}
\end{center}
\caption{Test 3. The domains $\Omega_k$ for $k=0$ (top, left),
intermediate domain during the line-search after $k=0$ ($\lambda=0.8^{15}$, top, middle),
$k=1$ (top, right), $k=2$ and $k=4$ (bottom).}
\label{fig:test7_Omega}
\end{figure}

We plotted some steps $\Omega_k$ in Figure \ref{fig:test7_Omega}
and we note that the topology changes during the iterations.

\begin{figure}[ht]
\begin{center}
  \includegraphics[width=9cm]{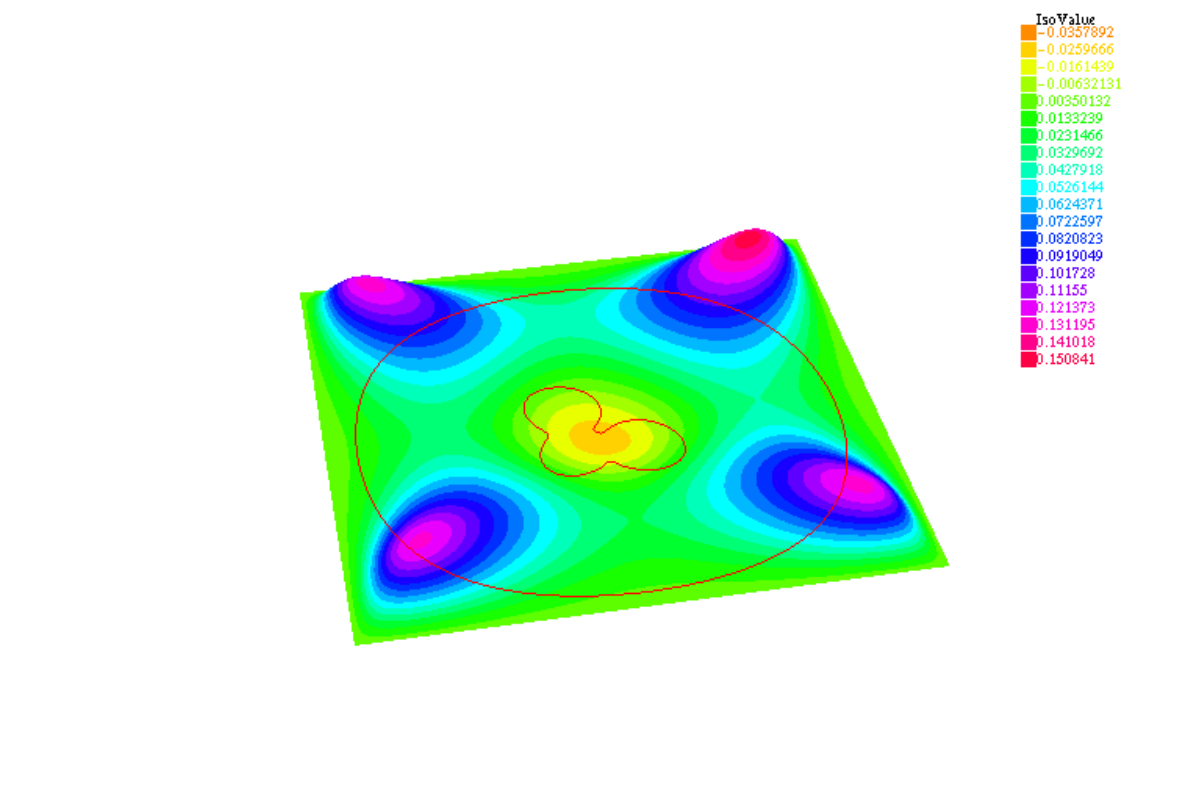}\\
  \includegraphics[width=6.5cm]{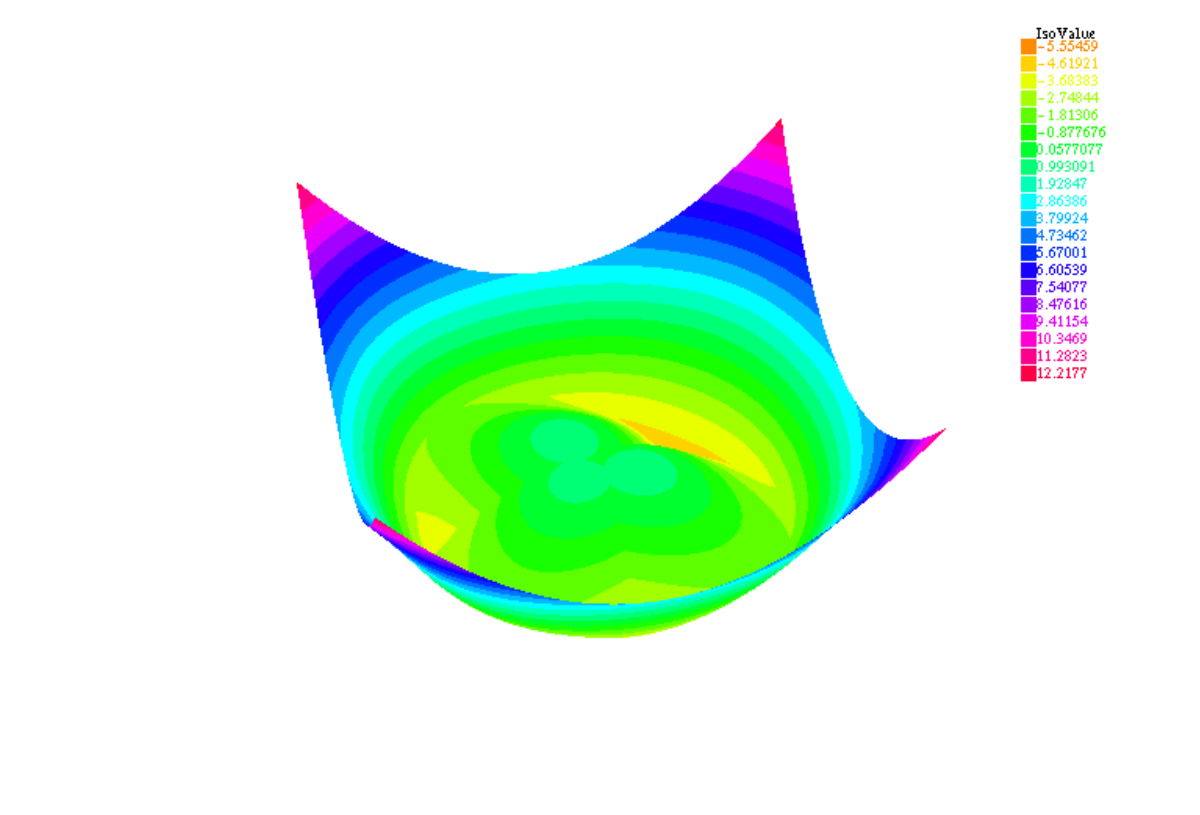}
  \includegraphics[width=6.5cm]{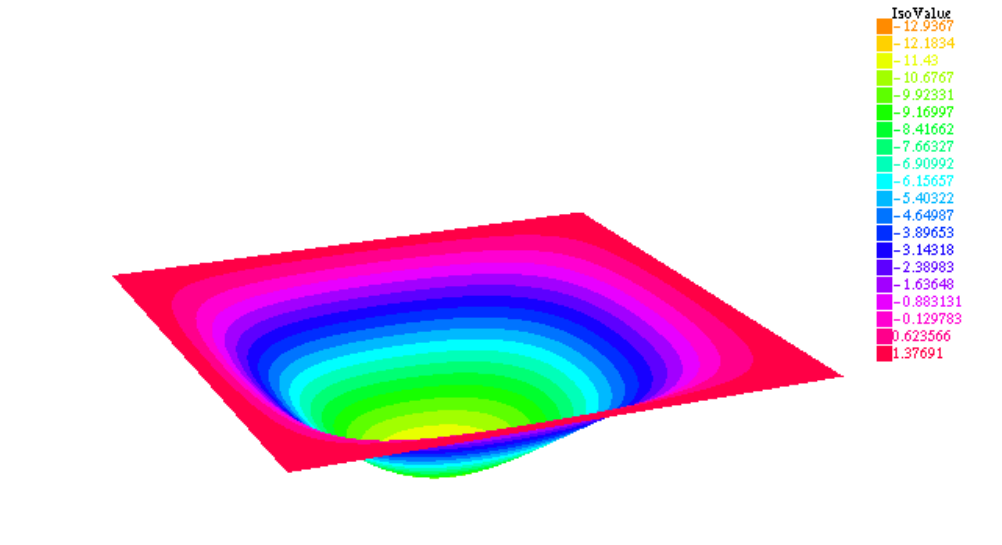}
\end{center}
\caption{Test 3. Final $y_h$ and $\partial\Omega_h$ (top).
  Final $g_h$ (left, bottom) and $u_h$ (right, bottom).}
\label{fig:test7_y}
\end{figure}

In Figure \ref{fig:test7_y}, we have plotted the final
$y_h,g_h,u_h:D\rightarrow\mathbb{R}$.
 The initial domain is as for Test 2, case b).
The value of the original cost at the initial iteration is a little bit different from
Test 2, case b) since the meshes are different.

\begin{figure}[ht]
\begin{center}
  \includegraphics[width=9cm]{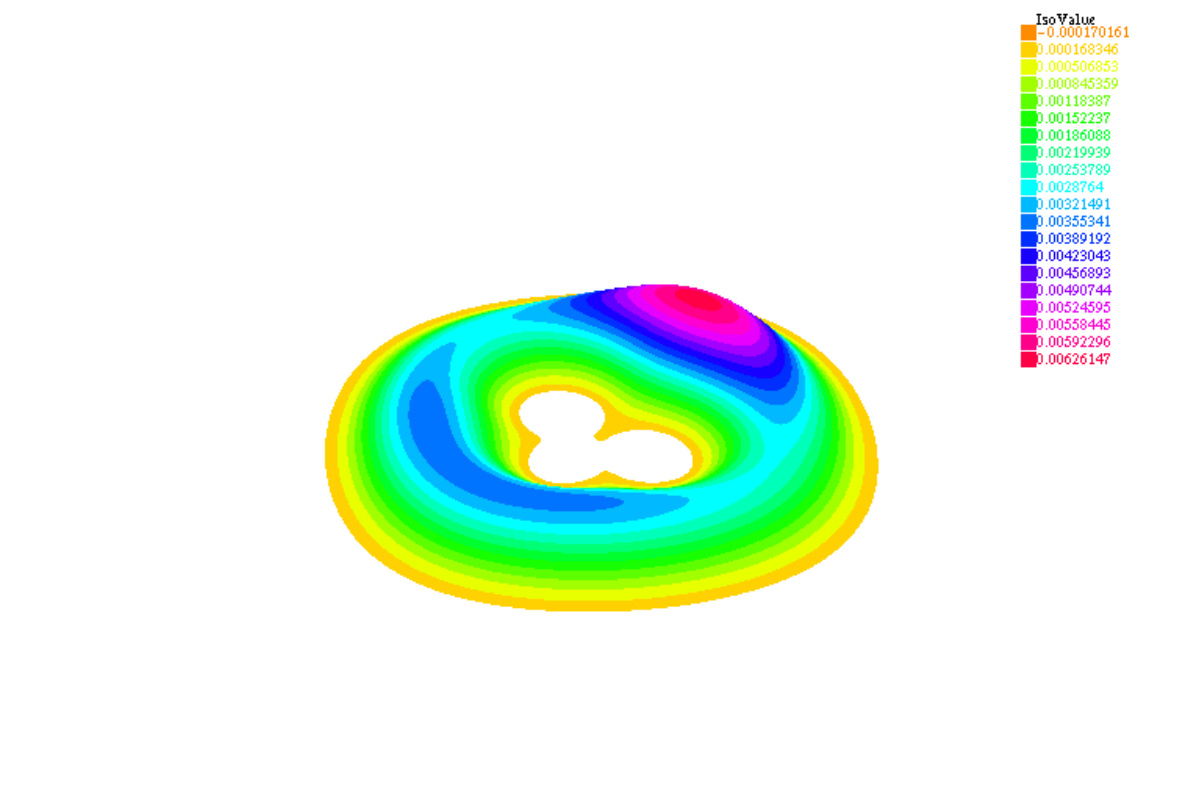}
\end{center}
\caption{Test 3. Solution of the clamped plate model in the final domain.}
\label{fig:test7_y1}
\end{figure}

\clearpage

\section*{Acknowledgement}
This work was partially supported by the French - Romanian cooperation
program ``ECO Math'', 2022-2023.

\end{document}